\pdfoutput=1
\documentclass[aos]{imsart}

\usepackage[T1]{fontenc}      
\usepackage[utf8]{inputenc}   
\usepackage{microtype}        

\usepackage{amsmath}
\usepackage{amsthm}
\usepackage{amssymb}
\usepackage{commath}
\usepackage{mathtools}
\usepackage{thmtools}
\usepackage{mathrsfs}

\usepackage[colorlinks,citecolor=blue,urlcolor=blue]{hyperref}
\usepackage{cleveref}
\usepackage[numbers]{natbib}
\usepackage[dvipsnames]{xcolor} 

\usepackage{comment}

\newcommand{\figpath}{}

\declaretheorem[Refname={Theorem,Theorems}]{theorem}
\numberwithin{theorem}{section}
\declaretheorem[style=definition,numberlike=theorem,Refname={Definition,Definitions}]{definition}

\declaretheorem[style=definition,numberlike=theorem,Refname={Remark,Remarks}]{remark}
\declaretheorem[numberlike=theorem,Refname={Lemma,Lemmas}]{lemma}
\declaretheorem[name=Corollary,numberlike=theorem,Refname={Corollary,Corollaries}]{corollary}
\declaretheorem[name=Proposition,numberlike=theorem,Refname={Proposition,Propositions}]{proposition}

\DeclarePairedDelimiterX\Set[2]{\lbrace}{\rbrace}%
{ #1 \,:\, #2 }                                         

\DeclareMathOperator*{\argmin}{arg\,min}

\newcommand{\R}{\mathbb{R}} 
\newcommand{\N}{\mathbb{N}} 
\newcommand{\Z}{\mathbb{Z}} 

\renewcommand{\b}[1]{\mathsf{#1}} 
\newcommand{\T}{\intercal} 
\DeclareMathOperator{\tr}{tr} 
\DeclareMathOperator{\diag}{diag} 

\newcommand{\E}{\mathbb{E}}     
\newcommand{\Var}{\mathrm{Var}}     
\newcommand{\GP}{\mathrm{GP}}   
\newcommand{\NN}{\mathrm{N}}   
\newcommand{\ML}{\textup{\textsf{ML}}}   
\newcommand{\CV}{\textup{\textsf{CV}}}   
\newcommand{\MSE}{\mathrm{MSE}} 
\newcommand{\Sob}{\mathrm{Sob}} 

\begin{document}

\begin{frontmatter}
        \title{Scale estimation and rate-unbiasedness for Gaussian processes under smoothness misspecification}
        \runtitle{Scale estimation under smoothness misspecification}
        
        \begin{aug}
        \author[A]{\fnms{Toni}~\snm{Karvonen}\ead[label=e1]{toni.karvonen@lut.fi}} and
        \author[B,C]{\fnms{Fran\c{c}ois}~\snm{Bachoc}\ead[label=e2]{francois.bachoc@univ-lille.fr}}
        \address[A]{School of Engineering Sciences, Lappeenranta--Lahti University of Technology LUT, Finland\printead[presep={ ,\ }]{e1}}
        \address[B]{Université de Lille, France\printead[presep={,\ }]{e2}}
        \address[C]{Institut universitaire de France (IUF), France}
        \end{aug}
        
        \begin{abstract}
                Gaussian process regression is used throughout statistics and machine learning for prediction and uncertainty quantification.
                A Gaussian process is specified by its mean and covariance functions. 
                Many covariance functions, including Matérns, have a smoothness parameter that is notoriously difficult to specify correctly or estimate from the data.
                In practice, the smoothness parameter is often selected more or less arbitrarily.
                We introduce rate-unbiasedness, a relaxed notion of asymptotic optimality which requires that the expected ratio of the mean-square error presumed by a potentially misspecified model and the true, but unknown, mean-square error remain bounded away from zero and infinity as more data are obtained.
                A rate-unbiased model provides uncertainty quantification that is of correct order of magnitude.                
                We then prove that scale estimation suffices for rate-unbiasedness in a variety of common settings.
                As estimation of the scale of a Gaussian process is routine and requires no optimisation, rate-unbiasedness can be achieved in many applications.
        \end{abstract}
        
        \begin{keyword}[class=MSC]
        \kwd[Primary ]{60G15}
        \kwd[; secondary ]{62G20, 46E22, 46E35, 65D15}
        \end{keyword}
        
        \begin{keyword}
        \kwd{Gaussian processes}
        \kwd{Matérn covariance}
        \kwd{parameter estimation}
        \kwd{misspecification}
        \end{keyword}
        
        \end{frontmatter}

\section{Introduction}

Gaussian processes constitute a model of choice in numerous fields within statistics and machine learning. 
For example, they have become popular for Bayesian modeling of complex computer experiments~\cite{Sacks1989,Santner2003} and are used in geostatistics under the name Kriging~\cite{cressie2015statistics, Stein1999}.
Gaussian processes have also been extensively applied to a variety of machine learning tasks~\cite{nickisch2008approximations,RasmussenWilliams2006}. In Bayesian nonparametric statistics, they exhibit favorable posterior contraction properties~\cite{VaartZanten2011,van2008rates,van2009adaptive}. 
Furthermore, Gaussian processes often serve as building blocks to create more complex structures, such as deep~\cite{bachoc2025posterior,castillo2024deep,damianou2013deep,finocchio2023posterior}, heavier-tailed~\cite{xu2017tukey}, and monotonic~\cite{Bachoc2019cMLE,DaVeiga2012GPineqconst,Pallavi2019BayesianShapeGPs,Wang2019DiffEq} processes.

In regression, a Gaussian process is observed at some locations of the input space, either exactly or with additive Gaussian noise. In this setting, the conditional mean function of the Gaussian process has a well-known explicit expression, and, even more, the entire conditional distribution of the process remains Gaussian~\cite[e.g.,][]{RasmussenWilliams2006,Stein1999}. Of course, the conditional mean and conditional covariance depend on the original mean function and covariance kernel of the Gaussian process. 
A central problem that has motivated a large amount of research at least since the 1980s~\cite{stein1988asymptotically} is to understand the properties of the conditional mean and covariance functions when the original mean and covariance used to construct them are \emph{misspecified}, meaning that they differ from those of the true underlying Gaussian process~\cite{bachoc2022asymptotically,stein1990uniform,stein1990bounds,Stein1999,zhang2015doesn}.

The equivalence of measures is a key concept in misspecification~\cite{ibragimov2012gaussian}.
Consider two Gaussian measures, one corresponding to the true pair of mean and covariance functions of the Gaussian process and the other to the misspecified pair used to compute the conditional mean and covariance. If these two measures are equivalent, then the computed conditional mean and covariance are asymptotically optimal in two senses \cite{stein1988asymptotically,stein1990bounds,Stein1999}: First, the ratio of mean-square prediction errors between the predictions obtained from the true and misspecified pairs tends to one as $n$, the number of observation locations, increases (here we always consider an increasing number of observation locations in a fixed input domain, a setting called fixed-domain asymptotics in \cite{Stein1999}).
Second, the ratio between the actual mean-square prediction error of the misspecified pair and the conditional variance computed from this misspecified pair (we call this the \emph{presumed} mean-square error), also tends to one. In short, under the equivalence of Gaussian measures prediction and uncertainty quantification remain asymptotically optimal even if the mean and covariance functions are misspecified.

Finding necessary and/or sufficient conditions for asymptotic equivalence of Gaussian measures is a classical problem in statistics and probability \cite{anderes2010consistent,ibragimov2012gaussian,skorokhod1973absolute}. 
The conditions that have been established are strong: essentially, the two pairs of mean and covariance functions must be very similar. 
For instance, consider the Matérn~\cite{genton2001classes,Stein1999} covariance
\begin{equation} \label{eq:matern-intro}
  K(x, y) = \sigma^2 \frac{2^{1-\nu}}{\Gamma(\nu)} \bigg( \frac{\sqrt{2\nu} \norm[0]{x - y}}{\lambda} \bigg)^{\!\nu} \mathcal{K}_{\nu} \bigg( \frac{\sqrt{2\nu} \norm[0]{x - y}}{\lambda} \bigg) \quad \text{ for } \quad x, y \in \R^d,
\end{equation}
where $\Gamma$ is the gamma function and $\mathcal{K}_\nu$ the modified Bessel function of the second kind.
The smoothness parameter $\nu > 0$ determines the regularity of the Gaussian random field, the correlation length $\lambda > 0$ controls how ``wiggly'' it is, and the scale parameter $\sigma > 0$ determines the magnitude of the variance.
Zero-mean Gaussian measures corresponding to Matérn covariances with parameters $(\nu_0, \lambda_0, \sigma_0)$ and $(\nu, \lambda, \sigma)$ are equivalent if and only if (a) $\nu_0 = \nu$ and (b) $\sigma_0^2 / \lambda_0^{2\nu} = \sigma^2 / \lambda^{2\nu}$ when $d \leq 3$ or $(\sigma_0, \lambda_0) = (\sigma, \lambda)$ when $d \geq 4$~\cite{anderes2010consistent,bolin2023equivalence,zhang2004inconsistent}.
That is, although theoretically important, settings in Gaussian process regression where fixed yet misspecified mean and covariance functions yield equivalent Gaussian measures are arguably not representative of statistics and machine learning practice. 
In practice, one typically first estimates the mean and covariance from the data and then uses these estimates to compute the conditional mean and covariance.
Only rarely does the estimation of the mean and covariance yield asymptotically optimal prediction and uncertainty quantification (see~\cite{putter2001effect} for one case).

In this article we introduce a relaxed notion of asymptotic optimality of uncertainty quantification and prove that it holds in great generality under covariance scale parameter estimation.
We assume that the data arise from a centered Gaussian process that is exactly observed at $n$ locations within a bounded domain $D \subset \R^d$.
The covariance of this data-generating process is a Sobolev kernel (\Cref{def:sobolev-kernel}) of smoothness $\nu_0$.
The class of Sobolev kernels extends that of Matérns.
However, the data-generating process is modelled as a centered Gaussian process with a Sobolev covariance kernel of smoothness $\nu \neq \nu_0$ and the conditional mean and covariance are computed based on this misspecified model.
The scale parameter $\sigma$ is estimated by maximum likelihood or cross-validation~\cite[Sec.\@~5.4]{RasmussenWilliams2006} and the mean-square error at $x \in D$ presumed by the model, $\MSE_n(x \mid \hat{\sigma}_n)$, is computed based on the misspecified covariance and an estimated scale parameter $\hat{\sigma}_n$.
All other covariance parameters, including smoothness, are fixed.
Prior results on scale parameter estimation have been obtained in~\cite{Karvonen2020, KorteStapff2024, Naslidnyk2024, SniekersVaart2015a, SniekersVaart2015b}.

\Cref{thm:rate-unbiased} shows that scale estimation suffices to make the presumed mean-square error \emph{rate-unbiased} if the observation locations are \emph{quasi-uniform} (\Cref{def:quasi-uniform}) and the model \emph{oversmooths} ($\nu \geq \nu_0$).
More precisely, for any $p \in (0, \infty)$ there are $c_1, c_2 > 0$ such that
\begin{equation} \label{eq:rate-unbiased-intro}
        c_1 \leq \frac{\mathbb{E}\big[\lVert \MSE_n(\cdot \mid \hat{\sigma}_n) \rVert_{L^p(D)}\big]}{\lVert \MSE_n^* \rVert_{L^p(D)}} \leq c_2
\end{equation}
for all $n$, where $\MSE_n^*(x)$ stands for the \emph{true} mean-square error, which is not computable without access to the true covariance function.
The expectation in~\eqref{eq:rate-unbiased-intro} is with respect to the data-generating Gaussian process and impacts the presumed mean-square error via the scale estimator.
Hence, in the simple practically realistic setting where one relies on the Mat\'ern model, fixes a (misspecified) smoothness parameter and estimates the scale parameter, the uncertainty quantification provided by the presumed mean-square error will have correct order of magnitude.
What makes this approach particularly convenient is that rate-unbiasedness is achieved \emph{without optimisation} since the maximum likelihood and cross-validation estimators of $\sigma$, the only parameter being estimated, are available in closed form [see~\eqref{eq:sigma-ML} and~\eqref{eq:cv-alternate}].
\Cref{thm:periodic} contains a limited extension of~\eqref{eq:rate-unbiased-intro} for periodic Sobolev covariances and $d = 1$ when the truth is at most twice as smooth as the model (i.e., the model can \emph{undersmooth}).
This theorem, which generalises results by Naslidnyk et al.\@~\cite{Naslidnyk2024} that apply to the Brownian motion, suggests that cross-validation should be preferred over maximum likelihood estimation under misspecification.
We note that parameter estimation for periodic covariances has been recently studied in~\cite{Chen2021,Petit2023}.
The proof of~\eqref{eq:rate-unbiased-intro} amounts to showing that $\E[\hat{\sigma}_n^2]$ blows up as $n^{2(\nu-\nu_0)/d}$ (recall that $\nu \geq \nu_0$).
\Cref{thm:estimator-asymptotics} contains this result of independent interest. 
In \Cref{sec:numerical}, we study numerically if $\nu_0$ can be estimated from an observed rate of increase of a scale estimator.

It would be preferable to obtain a stronger almost sure version of~\eqref{eq:rate-unbiased-intro}.
Unfortunately, we do not know how to do this for standard kernels and domains (if $D$ were a closed Riemannian manifold and the covariances Whittle--Matérns, the approach in~\cite{KorteStapff2024} would likely work).
Expectations of scale estimators are sums of ratios of mean-square errors, which can be interpreted as \emph{worst-case errors in the reproducing kernel Hilbert spaces} of the true and presumed covariances [see~\eqref{eq:wce-mse-1},~\eqref{eq:wce-mse-2},~\eqref{eq:sigma-ML-E-app}, and~\eqref{eq:sigma-CV-E-app}].
Tools and techniques from kernel-based approximation~\cite{Wendland2005, Iske2018} and information-based complexity~\cite{Novak1988,NovakWozniakowski2008,NovakWozniakowski2010} can be then exploited to understand the asymptotic behaviour of these worst-case errors.
In contrast, an almost sure result would require bounding approximation errors for realisations from a Gaussian process, which is much more difficult than bounding worst-case errors (lower bounds are particularly challenging).
The assumption that the observation locations are quasi-uniform is relatively weak as it does not impose a particular design or require any type of stratification.

While this article focuses on noise-free data stemming from a Gaussian process, we point out that there is substantial recent Bayesian nonparametrics literature that considers closely related settings~\cite{castillo2024deep,ghosal2017fundamentals,HadjiSzabo2021,szabo2025adaptation,szabo2015frequentist,van2008rates,van2009adaptive}. 
In these references, Gaussian process models and Gaussian sequence models are shown to have various beneficial properties, such as adaptation for functional estimation or reliability of credible sets, when their parameters are estimated or given additional priors.
This is similar to the conclusion of this article that estimating the scale parameter from the data allows reliable uncertainty quantification. 
In Bayesian nonparametrics the data are typically assumed to be noisy and to stem from an unknown but fixed function or sequence.
In this article, the fixed function is replaced by a Gaussian process and, given a realisation from this process, the data are deterministic.
Consequently, the mathematical techniques that we use, which rely on kernel-based approximation and information-based complexity, differ greatly from those used in the above references.

\subsection{Notation and conventions}

For non-negative sequences $(a_n)_{n=1}^\infty$ and $(b_n)_{n=1}^\infty$ we write $a_n \lesssim b_n$ if there is $C \geq 0$ such that $a_n \leq C \, b_n$ for all $n \geq 1$.
We write $a_n \asymp b_n$ if there is $C \geq 1$ such that $C^{-1} b_n \leq a_n \leq C \, b_n$.
These notations are equivalent to $a_n = O(b_n)$ and $a_n = \Theta(b_n)$.
When $a_n$ and $b_n$ are non-negative random variables, we write $a_n \asymp_{\mathbb{P}} b_n$ if
\begin{equation} \label{eq:asymp-in-prob}
        \limsup_{n \to \infty} \mathbb{P} ( a_n \le \varepsilon \, b_n ) \to 0 
        \quad \text{ and } \quad 
        \limsup_{n \to \infty} \mathbb{P} ( a_n \ge \varepsilon^{-1} b_n )
        \to
        0 \quad \text{ as } \quad \varepsilon \to 0 .
\end{equation}
We consider positive-definite kernels on a domain $D \subseteq \R^d$.
We say that $K \colon D \times D \to \R$ is a positive-definite kernel on $D$ if $K$ is symmetric, which is to say that $K(x, y) = K(y ,x)$ for all $x, y \in D$, and if the covariance matrix $\b{K}_n = (K(x_i, x_j))_{i,j=1}^n$ is strictly positive-definite for all $n \in \N$ and all pairwise distinct locations $x_1, \ldots, x_n \in D$.
This implies that the covariance matrix is invertible.
A kernel is positive-semidefinite if the matrix $\b{K}_n$ is positive-semidefinite.

\section{Setting and preliminaries} \label{sec:setting}

This section describes the Gaussian process modelling setting we consider and reviews necessary preliminaries on Sobolev spaces and related topics.

\subsection{Gaussian process modelling} \label{sec:GP-modelling}

\begin{figure}
        \centering
        \includegraphics[width=\textwidth]{\figpath 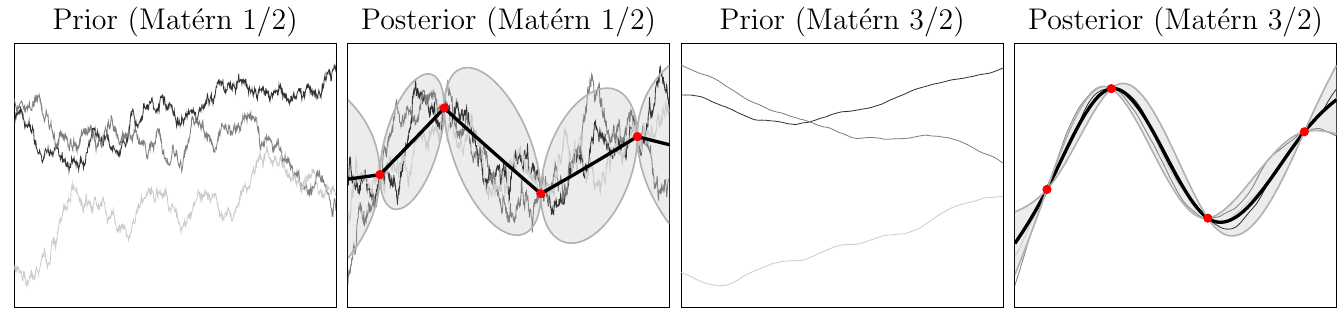}
        \caption{Samples from univariate Gaussian process priors and posteriors defined by the Matérn kernel in~\eqref{eq:matern-intro} with $\nu \in \{1/2, 3/2\}$ and $\lambda = 1$. The red dots are the data points $(x_i, y_i)$, the black line is the posterior mean function, and the shaded region represents the $95\%$ credible intervals around the mean.}
        \label{fig:GP-illustration}
\end{figure}

Standard references on Gaussian processes include \citep{Gramacy2020, RasmussenWilliams2006, Santner2003}.
Suppose that we have observations $y_1, \ldots, y_n \in \R$ corresponding to pairwise distinct locations $x_1, \ldots, x_n \in D$.
To predict observations at unseen locations, we can postulate that the observations arise from a realisation of a random process $X$ and apply Bayesian methodology to infer the most likely observations given the data $\mathcal{D}_n = \{(x_i, y_i)\}_{i=1}^n$.
Choosing a zero-mean Gaussian process prior $X \sim \GP(0, K)$ with a positive-definite covariance $K \colon D \times D \to \R$ permits closed-form conditioning.
We defer detailed discussion on covariance kernels to \Cref{sec:sobolev}.
For notational simplicity we use a zero-mean prior throughout this article; it would be straightforward to relax this assumption.
Under this prior the vector $\mathsf{y}_n = (y_1, \ldots, y_n) \in \R^n$ that collects observations is a zero-mean Gaussian random vector with positive-definite covariance matrix $\b{K}_n = (K(x_i, x_j))_{i,j=1}^n$.
The equations for Gaussian conditioning yield the Gaussian posterior process $X \mid \mathcal{D}_n \sim \GP(\mu_n, C_n)$ whose mean and covariance are given by
\begin{equation} \label{eq:GP-mean-cov}
        \mu_n(x) = \b{k}_n(x)^\T \b{K}_n^{-1} \b{y}_n \quad \text{ and } \quad C_n(x, y) = K(x, y) - \b{k}_n(x)^\T \b{K}_n^{-1} \b{k}_n(y),
\end{equation}
where $\b{k}_n(x) = (K(x, x_1), \ldots, K(x, x_n)) \in \R^n$.
\Cref{fig:GP-illustration} shows two Gaussian process priors and the resulting posteriors.
We shall mostly work with the conditional variance
\begin{equation} \label{eq:posterior-var}
        V_n(x) = C_n(x, x) = K(x, x) - \b{k}_n(x)^\T \b{K}_n^{-1} \b{k}_n(x) .
\end{equation}
At each $x \in D$, the variance equals the mean-square error as \emph{presumed by the model}:
\begin{equation} \label{eq:MSE}
        \MSE_n(x) = \E_X[ X(x) - \mu_n(x) ]^2 = V_n(x),
\end{equation}
where $\E_X$ indicates that the expectation is taken under the assumption that the observations come from the Gaussian process $X \sim \GP(0, K)$, which is to say that $\b{y}_n \sim \NN(0, \b{K}_n)$.

Rarely, if ever, are the observations generated by the Gaussian process $X$, or a Gaussian process at all.
Suppose that in reality the observations arise from an unknown Gaussian process $X_0 \sim \GP(0, K_0)$ whose covariance $K_0$ need not equal $K$.
Then $\b{y}_n \sim \NN(0, \b{K}_{0,n})$, where $\b{K}_{0,n} = (K_0(x_i, x_j))_{i,j=1}^n$ is the covariance matrix for the true data-generating process.
It follows that the \emph{true mean-square error} of the conditional mean $\mu_n(x)$ in~\eqref{eq:GP-mean-cov} is
\begin{equation} \label{eq:MSE-true}
\begin{split}
        \MSE_n^*(x) &= \E[ X_0(x) - \mu_n(x) ]^2 \\
        &= \E\big[ X_0(x)^2 - 2 X_0(x) \mu_n(x) + \mu_n(x)^2 \big] \\        
        &= K_0(x, x) - 2 \b{k}_{0,n} \b{K}_n^{-1} \b{k}_n(x) + \b{k}_n(x)^\T \b{K}_n^{-1} \b{K}_{0,n} \b{K}_n^{-1} \b{k}_n(x) ,
\end{split}
\end{equation}
where we used $\mu_n(x) = \b{k}_n(x)^\T \b{K}_n^{-1} \b{y}_n$ and the covariances $\E[X_0(x) \b{y}_n] = \b{k}_{0,n}(x)$ and \smash{$\E[\b{y}_n \b{y}_n^\T] = \b{K}_{0,n}$}.
However, because the true covariance $K_0$ is unknown, one cannot access \smash{$\MSE_n^*(x)$} and must instead substitute $\MSE_n(x)$ in~\eqref{eq:MSE} as a measure of error and uncertainty.
The mean-square errors may differ greatly in how they behave, and to use $\MSE_n(x)$ may be accordingly highly misleading.
The following theorem illustrates what we mean.
The theorem uses the concept of a Sobolev kernel that we shall introduce in \Cref{def:sobolev-kernel}.
A Matérn kernel of order $\nu$ in~\eqref{eq:matern-intro} is a Sobolev kernel of order $\nu$. 
See \Cref{sec:proofs-setting} for a proof of \Cref{thm:MSE-rates-informal}.

\begin{theorem} \label{thm:MSE-rates-informal}
Let $p \in (0, \infty)$ and $\nu \geq \nu_0 > 0$.
Suppose that $D$ is a bounded open convex set and $K$ and $K_0$ are Sobolev kernels of orders $\nu$ and $\nu_0$.
If the sequence $(x_i)_{i=1}^\infty \subset D$ is quasi-uniform, then
\begin{equation} \label{eq:MSE-rates}
        \lVert \MSE_n \rVert_{L^p(D)} \asymp n^{-2\nu/d} \quad \text{ and } \quad \lVert \MSE_n^* \lVert_{L^p(D)} (x) \asymp n^{-2\nu_0/d} .
\end{equation}
\end{theorem}

\Cref{thm:MSE-rates-informal} states that a model that oversmooths the truth (i.e., $\nu > \nu_0$) is liable to overconfidence as its presumed mean-square error tends to zero much faster than the true error.
Our goal is to show that equipping the model with a scale parameter and using maximum likelihood estimation or cross-validation to select this parameter eliminates the problem.

\begin{remark}
In \Cref{thm:MSE-rates-informal} and other results of this article we assume that $D$ is convex.
We use this assumption due to its simplicity rather than its necessity.
Convexity can be replaced with the less restrictive interior cone condition and an assumption that the domain have Lipschitz boundary, both of which are assumptions that commonly appear in the theory of Sobolev spaces.
The more general assumptions are used in, for example, \cite{Karvonen2020, WangJing2022}.
\end{remark}

\subsection{Parameter estimation and rate-unbiasedness}

In practice, the Gaussian process model is usually specified in terms of a covariance kernel $K_\theta$ with some parameters $\theta \in \Theta$ that are estimated from the data.
The scale, correlation length, and smoothness parameters $\sigma$, $\lambda$, and $\nu$ of the Matérn model~\eqref{eq:matern-intro} are among the most frequently estimated kernel parameters.
The kernel with estimated parameters is subsequently plugged in the posterior equations~\eqref{eq:GP-mean-cov}.
Maximum likelihood estimation and cross-validation, the two most popular approaches to plug-in parameter estimation, require troublesome and computationally expensive optimisation.

The Gaussian process model is said to be \emph{well-specified} if the true data-generating process $X_0 \sim \GP(0, K_0)$ falls within the parametrisation, in that $K_0 \in \{K_\theta\}_{\theta \in \Theta}$.
In this case it is possible (though not necessarily easy) to recover the true covariance kernel.
However, some kernel parameters are difficult to estimate and thus tend to be fixed beforehand, a problem particularly associated with the Matérn smoothness parameter $\nu$ (see~\cite{Geoga2023} for recent work on the topic).
Even if all parameters are estimated, the parametric model may be too inflexible or constrained to contain $K_0$, or the true process may be so esoteric that it is not contained in any standard family of processes.
For example, for every $\lambda > 0$ the Gaussian covariance kernel
\begin{equation} \label{eq:gaussian-kernel}
        K(x, y) = \sigma^2 \exp\bigg(\!-\frac{\lVert x - y \rVert^2}{2\lambda^2} \bigg), 
\end{equation}
that remains popular in machine learning, induces a Gaussian process with paths so smooth that it is unrealistic to expect any physically relevant process to come from this parametric family~\citep[p.\@~55]{Stein1999}.
One should therefore always hold it more likely than not that the model is \emph{misspecified}, which is to say that the true kernel is not an element of the parametric class: $K_0 \notin \{K_\theta\}_{\theta \in \Theta}$ \cite{Bachoc2018}.
Although the model cannot fully capture the behaviour of the true data-generating process in the misspecified setting, one does not have to abandon all hope.
Let
\begin{equation} \label{eq:MSE-implied-theta}
        \MSE_n(x \mid \theta) = V_n(x \mid \theta) = K_\theta(x, x) - \b{k}_{\theta,n}^\T \b{K}_{\theta,n}^{-1} \b{k}_{\theta,n}(x) 
\end{equation}
denote the presumed mean-square error given kernel parameters $\theta$.
Here $\b{k}_{\theta,n}(x)$ and $\b{K}_{\theta,n}$ are defined as $\b{k}_n(x)$ and $\b{K}_n$ in \Cref{sec:GP-modelling} but with $K = K_\theta$.
We introduce an extension of the notion of unbiasedness that relaxes asymptotic optimality.

\begin{definition}[Rate-unbiasedness]
Let $(x_i)_{i=1}^\infty$ be a sequence of observation locations and \smash{$\hat{\theta}_n$} a kernel parameter estimator based on the data $\mathcal{D}_n$. We say that the presumed mean-square error \smash{$\MSE_n(\cdot \mid \hat{\theta}_n)$} in~\eqref{eq:MSE-implied-theta} is \emph{rate-unbiased} with $p \in (0, \infty]$ if its expected $L^p(D)$-norm behaves asymptotically as the true mean-square error $\MSE_n^*$ in~\eqref{eq:MSE-true}:
\begin{equation} \label{eq:rate-unbiased}   
                    \E\big[ \, \lVert \MSE_n(\cdot \mid \hat{\theta}_n) \rVert_{L^p(D)} \, \big] \asymp \lVert \MSE_n^* \rVert_{L^p(D)} . 
\end{equation}  
If the presumed mean-square error is rate-unbiased for no $p \in (0, \infty]$, then it is \emph{rate-biased}.
\end{definition}

In \Cref{thm:MSE-rates-informal} we saw that, in the case of a Matérn model with misspecified smoothness and the trivial estimator $\hat{\theta}_n = (\sigma, \lambda, \nu)$ for all $n \geq 1$, the presumed mean-square error is rate-biased if $\nu > \nu_0$.
Below we give a more interesting example.
See \Cref{sec:proofs-setting} for a proof.

\begin{proposition} \label{prop:rate-unbiased-lambda}
Suppose that $D$ is a bounded open convex set.
Let $K$ be a Matérn kernel in~\eqref{eq:matern-intro} with fixed $\sigma$ and $\nu$ and $K_0$ a Matérn kernel with parameters $\sigma_0$, $\lambda_0$, and $\nu_0$.
Suppose that $(x_i)_{i=1}^\infty$ is any sequence of observation locations.
If $\hat{\lambda}_n$ is an estimator of $\lambda$ such that $0 < C_1 \leq \hat{\lambda}_n \leq C_2 < \infty$ for all $n \geq 1$, then the presumed mean-square error is rate-unbiased with any $p \in (0, \infty)$ if $\nu_0 = \nu$ and rate-biased if $\nu > \nu_0$.
\end{proposition}

The proposition states that an estimator of the correlation length that is bounded away from zero and infinity does not affect rate-unbiasedness: If smoothness is correctly specified, any such estimator induces rate-unbiasedness; if the model oversmooths, the presumed mean-square error is rate-biased.
To ensure that Gaussian process modelling performs well when $n$ is small it is essential to use a sensible correlation length estimator.
\Cref{prop:rate-unbiased-lambda} thus tells us that rate-unbiasedness is a rather rather weak, yet important, property.
Without rate-unbiasedness uncertainty quantification is guaranteed to be unreliable (recall \Cref{thm:MSE-rates-informal}) but rate-unbiasedness alone is not enough for a model to perform well: practically important parameters, such as the correlation length, can be misspecified and the constants in~\eqref{eq:rate-unbiased} that are hidden by asymptotic notation need not be close to one.
We refer to~\cite{Stein1990} and~\cite[Ch.\@~3]{Stein1999} for more discussion on the relation between the true and presumed mean-square errors.

\subsection{Scale parameter estimation} \label{sec:scale-estimation}

Estimation of a \emph{scale parameter} is an important special case.
Let $X \sim \GP(0, K)$ be a fixed Gaussian process and suppose that the process
\begin{equation} \label{eq:scale-parametrisation}
        X_\sigma = \sigma X \sim \GP(0, \sigma^2 K)
\end{equation}
used to model the observations has a positive scale parameter $\sigma$.
Under this model the mean and variance of the posterior process $X_\sigma \mid \mathcal{D}_n$ are (note that the mean is independent of $\sigma$)
\begin{equation*}
        \mu_n(x) = \b{k}_n(x)^\T \b{K}_n^{-1} \b{y}_n \quad \text{ and } \quad V_n^\sigma(x) = \sigma^2 V_n(x) = \sigma^2 \big[ K(x, x) - \b{k}_n(x)^\T \b{K}_n^{-1} \b{k}_n(x) \big] .
\end{equation*}
The mean square-error implied by the model is now
\begin{equation} \label{eq:MSE-implied-sigma}
\MSE_n(x \mid \sigma) = \sigma^2 \MSE_n(x) = \sigma^2 V_n(x) .
\end{equation}
It is convenient to use maximum likelihood or cross-validation to set the scale as the resulting estimators 
are available in closed form in terms of conditional means and variances~\citep[Sec.\@~3]{Bachoc2013}.
The maximum likelihood estimator for $\sigma^2$ is
\begin{equation} \label{eq:sigma-ML}
        \hat{\sigma}_{\ML,n}^2 = \frac{1}{n} \b{y}_n^\T \b{K}_n^{-1} \b{y}_n = \frac{1}{n} \sum_{k=1}^n \frac{( y_k - \mu_{k-1}(x_k))^2}{V_{k-1}(x_k)} ,
\end{equation}
where the expansion, which is written in terms of conditional means and variances given datasets consisting of the $k = 0, \ldots, n-1$ first data points, is relatively well-known; see \cite[Sec.\@~4.2.2]{XuStein2017} or \cite[Prop.\@~7.5]{KarvonenOates2023}.
Here $\mu_0 \equiv 0$ and $V_0(x) = K(x, x)$.
The expansion can be proved by applying the block matrix inversion formula and using the linear-algebraic expression for the conditional variance in~\eqref{eq:posterior-var}.
The leave-one-out cross-validation estimator is
\begin{equation} \label{eq:sigma-CV}
        \sigma_{\CV,n}^2 = \frac{1}{n} \sum_{k=1}^n \frac{( y_k - \mu_{n \setminus k}(x_k))^2}{V_{n \setminus k}(x_k)} ,
\end{equation}
where the subscript indicates that the conditional mean and variance are formed using the dataset $\mathcal{D}_n \setminus \{(x_k, y_k)\}$ that contains all but the $k$th data point.
Note the similarity of the two estimators.
It is not a difficult exercise in linear algebra~\citep[Sec.\@~5.4.2]{RasmussenWilliams2006} to show that the cross-validation estimator admits the computationally useful alternate expression
\begin{equation} \label{eq:cv-alternate}
        \sigma_{\CV,n}^2 = \frac{1}{n} \b{y}_n^\T \b{K}_n^{-1} (\diag \b{K}_n^{-1})^{-1} \b{K}_n^{-1} \b{y}_n ,
\end{equation}
where $\diag \b{A}$ stands for the diagonal matrix whose diagonal coincides with that of $\b{A}$.

The model~\eqref{eq:scale-parametrisation} parametrised by scale is well-specified if $K_0 = \sigma_0^2 K$ for some $\sigma_0^2$.
In this case, both the maximum likelihood and cross-validation estimators of $\sigma^2$ are unbiased because
\begin{equation} \label{eq:sigma-ML-unbiased}
        \E[ \hat{\sigma}_{\ML,n}^2 ] = \frac{1}{n} \tr( \b{K}_n^{-1} \E[ \b{y}_n \b{y}_n^\T ] ) = \frac{1}{n} \tr( \b{K}_n^{-1} \b{K}_{0,n} ) = \frac{1}{n} \tr(\sigma_0^2 \b{Id}_n) = \sigma_0^2
\end{equation}
and
\begin{equation*}
        \E[ \hat{\sigma}_{\CV,n}^2 ] = \frac{1}{n} \sum_{k=1}^n \frac{\E[ y_k - \mu_{n \setminus k}(x_k)]^2}{V_{n \setminus k}(x_k)} = \frac{1}{n} \sum_{k=1}^n \frac{\sigma_0^2 V_{n \setminus k}(x_k)}{V_{n \setminus k}(x_k)} = \sigma_0^2 . 
\end{equation*}

Our goal is to show that the scale estimators $\hat{\sigma}_{\ML,n}$ and $\hat{\sigma}_{\CV,n}$ induce rate-unbiasedness in the setting of \Cref{thm:MSE-rates-informal}.
By~\eqref{eq:MSE-implied-sigma},
\begin{equation*}
        \E\big[ \, \lVert \MSE_n(\cdot \mid \hat{\sigma}_n) \rVert_{L^p(D)} \, \big] = \E[ \hat{\sigma}_n^2 ] \cdot \lVert \MSE_n \rVert_{L^p(D)} 
\end{equation*}
for any scale estimator $\hat{\sigma}_n$ and any $p \in (0, \infty]$. 
It follows from~\eqref{eq:MSE-rates} that 
\begin{equation} \label{eq:scale-est-required-rate}
        \E[ \hat{\sigma}_n^2 ] \asymp \frac{\lVert \MSE_n^* \rVert_{L^p(D)}}{ \lVert \MSE_n \rVert_{L^p(D)}} \asymp \frac{n^{-2\nu_0/d}}{n^{-2\nu/d}} = n^{2(\nu-\nu_0)/d}
\end{equation}
is a sufficient condition for rate-unbiasedness when $\nu \geq \nu_0$.
Our proof of rate-unbiasedness therefore consists of establishing this asymptotic for the two scale estimators. 
Note that the right-hand of~\eqref{eq:scale-est-required-rate} is of constant order or blows up since we assume $\nu \geq \nu_0$.
Next we properly introduce the notions that have already been used in \Cref{thm:MSE-rates-informal}.

\subsection{Sobolev spaces and kernels} \label{sec:sobolev}

Each positive-semidefinite kernel $K \colon D \times D \to \R$ induces a unique \emph{reproducing kernel Hilbert space} (RKHS) $H(K)$, a Hilbert space of functions $f \colon D \to \R$ with an inner product $\langle \cdot , \cdot \rangle_{H(K)}$ such that $K(\cdot, x) \in H(K)$ for every $x \in D$ and 
\begin{equation} \label{eq:reproducing-property}
        f(x) = \langle f, K(\cdot, x) \rangle_{H(K)} \quad \text{ for all } \quad f \in H(K) \text{ and } x \in D. 
\end{equation}
Equation~\eqref{eq:reproducing-property} is called the \emph{reproducing property}.
See \citep{Berlinet2004, Paulsen2016} for an introduction to RKHSs.

The RKHS of a stationary kernel is determined by its spectral density.
If $K(x, y) = \Phi(x - y)$ is a stationary kernel defined by an integrable and continuous $\Phi \colon \R^d \to \R$, then the RKHS of $K$ on $\R^d$, $H(K, \R^d)$, consists of those square-integrable functions $f$ whose Fourier transforms $(\mathcal{F}f)(\omega) = (2\pi)^{-d/2} \int_{\R^d} f(x) \exp(-\mathrm{i} x^\T \omega ) \dif x$ satisfy
\begin{equation} \label{eq:RKHS-fourier}
        \lVert f \rVert_{H(K, \R^d)}^2 = \frac{1}{(2\pi)^{d/2}} \int_{\R^d} \frac{\lvert (\mathcal{F}f)(\omega) \rvert^2}{(\mathcal{F} \Phi)(\omega)} \dif \omega < \infty .
\end{equation}
See, for example, Theorem~10.12 in \citep{Wendland2005}.
The Fourier transform $\mathcal{F} \Phi$, which is non-negative by Bochner's theorem, is called the spectral density of $K$.
On a proper subset $D$ of $\R^d$, the RKHS contains those functions which admit an extension onto $\R^d$ that satisfies~\eqref{eq:RKHS-fourier}.
That is, 
\begin{equation} \label{eq:RKHS-restricted}
        H(K) = \Set{f \colon D \to \R}{\text{ $f = f_e|_D$ for some $f_e \colon \R^d \to \R$ that satisfies~\eqref{eq:RKHS-fourier}}} .
\end{equation}
The norm $\lVert f \rVert_{H(K)}$ equals the minimum over the $H(K,\R^d)$-norms in~\eqref{eq:RKHS-fourier} of all possible extensions~\cite[Cor.\@~5.8]{Paulsen2016}.
The spectral density of the Matérn covariance in~\eqref{eq:matern-intro} of order $\nu$ is
\begin{equation} \label{eq:matern-FT}
        (\mathcal{F} \Phi)(\omega) = \sigma^2 \frac{2^{d/2} \Gamma(\nu + d/2)}{\Gamma(\nu)} \bigg(\frac{2\nu}{\lambda^2} \bigg)^\nu \bigg( \frac{2\nu}{\lambda^2} + \lVert \omega \rVert^2 \bigg)^{-(\nu+d/2)} .
\end{equation}
See, for example, page~49 in \citep{Stein1999} and bear the different parametrisations.

The Sobolev space $H^\alpha(\R^d)$ of order $\alpha > 0$ is a Hilbert space that consists of those functions $f \in L^2(\R^d)$ that satisfy
\begin{equation} \label{eq:Bessel-norm}
                 \lVert f \rVert_{H^\alpha(\R^d)}^2 = \frac{1}{(2\pi)^{d/2}}\int_{\R^d} \lvert (\mathcal{F} f)(\omega) \rvert^2 (1 + \lVert \omega \rVert^2)^\alpha \dif \omega < \infty.
\end{equation}
These spaces are conventionally called Bessel potential spaces but we eschew this terminology for simplicity.
On an arbitrary subset $D$ of $\R^d$, the space $H^\alpha(D)$ is defined analogously to the RKHS in~\eqref{eq:RKHS-restricted} as the space of functions which have extensions that satisfy~\eqref{eq:Bessel-norm}.
The standard definition of a Sobolev space of integer order, $W^{\alpha,2}(D)$, on a measurable set $D$ is via weak derivatives $\mathrm{D}^{\b{k}} f = \mathrm{D}^{k_1}_1 \cdots \mathrm{D}^{k_d}_d f$ and the norm
\begin{equation} \label{eq:sobolev-norm}
        \lVert f \rVert_{W^{\alpha,2}(D)}^2 = \sum_{ \lvert \b{k} \rvert \leq \alpha} \lVert \mathrm{D}^{\b{k}} f \rVert_{L^2(D)}^2 ,
\end{equation}
where $\b{k} \in \N_0^d$ are non-negative multi-indices.
Normed spaces $H$ and $F$ are \emph{norm-equivalent} if they are equal as sets and if there are non-negative constants $C_1$ and $C_2$ such that
\begin{equation*}
        C_1 \lVert f \rVert_F \leq \lVert f \rVert_H \leq C_2 \lVert f \rVert_F \quad \text{ for all } \quad f \in H.
\end{equation*}
If $D$ is sufficiently regular (e.g., convex) and $\alpha$ is integer, the Sobolev spaces in~\eqref{eq:Bessel-norm} and~\eqref{eq:sobolev-norm} are norm-equivalent~\citep[Cor.\@~10.48]{Wendland2005}.
On $D = \R^d$, this follows from the binomial theorem and the properties of the Fourier transform.
For example, in the one-dimensional case
\begin{align*}
        \lVert f \rVert_{H^\alpha(\R)}^2 = \frac{1}{\sqrt{2\pi}}\int_{\R} \lvert (\mathcal{F} f)(\omega) \rvert^2 (1 + \omega^2)^\alpha \dif \omega &= \sum_{k=0}^\alpha \binom{\alpha}{k} \frac{1}{\sqrt{2\pi}} \int_\R \lvert (\mathcal{F}f)(\omega) \rvert^2 \omega^{2k} \dif \omega \\
&= \sum_{k=0}^\alpha \binom{\alpha}{k} \lVert \mathrm{D}^k f \rVert_{L^2(\R)}^2
\end{align*}
shows that the norms $\lVert \cdot \rVert_{H^\alpha(\R)}$ and $\lVert \cdot \rVert_{W^{\alpha,2}(\R)}$ are equivalent.
In the last equality we used $(\mathcal{F} \, \mathrm{D}^k f)(\omega) = (\mathrm{i} \omega)^{k} (\mathcal{F}f)(\omega)$ and the Plancherel theorem.

Inserting the spectral density of a Matérn in~\eqref{eq:matern-FT} to the Fourier characterisation~\eqref{eq:RKHS-fourier} of the RKHS of a stationary kernel shows that the RKHS of a Matérn of order $\nu$ is norm-equivalent to $H^\alpha(D)$ for $\alpha = \nu + d/2$ on any $D \subseteq \R^d$.
Our results apply to Sobolev kernels that generalise the Matérn class in~\eqref{eq:matern-intro} in the sense that their RKHSs are norm-equivalent to Sobolev spaces.

\begin{definition}[Sobolev kernel] \label{def:sobolev-kernel}
A positive-definite kernel $K \colon D \times D \to \R$ on a set $D \subseteq \R^d$ is a \emph{Sobolev kernel} of order $\nu > 0$ on $D$ if the RKHS of $K$ on $D$ is norm-equivalent to the Sobolev space $H^\alpha(D)$ with $\alpha = \nu + d/2$.
In this case we write $K \in \Sob(\nu)$.
\end{definition}

\begin{figure}
        \centering
        \includegraphics[width=\textwidth]{\figpath 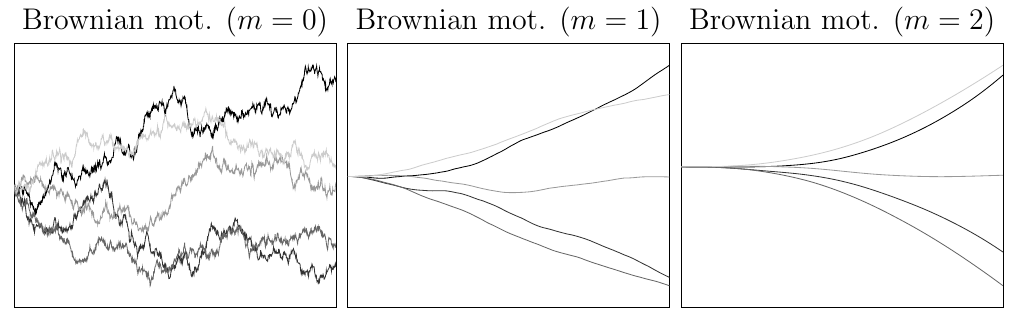}
        \caption{Samples from the Brownian motion ($m=0$), the integrated Brownian motion ($m=1$), and the twice integrated Brownian motion ($m=2$). Covariance kernels of these Gaussian processes are given in~\eqref{eq:ibm-kernel}.}
        \label{fig:bm-samples}
\end{figure}

Note that a Sobolev kernel of order $\nu$ is often~\citep[e.g.,][]{Karvonen2020} defined as a kernel whose RKHS is norm-equivalent to \smash{$H^\nu(D)$}, rather than $H^{\smash{\nu + d/2}}(D)$.
For our purposes the present definition is more convenient.
By~\eqref{eq:RKHS-fourier} and~\eqref{eq:Bessel-norm}, any stationary kernel whose spectral density satisfies
\begin{equation*}
        C_1 (1 + \lVert \omega \rVert^2)^{-(\alpha+d/2)} \leq (\mathcal{F} \Phi)(\omega) \leq C_2 (1 + \lVert \omega \rVert^2)^{-(\alpha+d/2)}
\end{equation*}
for some positive $C_1$ and $C_2$ and all $\omega \in \R^d$ is a Sobolev kernel of order $\alpha$.
However, a Sobolev kernel need not be stationary.
For example, the non-stationary released Brownian motion kernel $K(x, y) = 1 + \min\{x, y\}$ is a Sobolev kernel of order $\nu = 1/2$ on the interval $D = [0, 1]$ because its RKHS has the norm $\smash{\lVert f \rVert_{H(K)}^2 = f(0)^2 + \int_0^1 [ \mathrm{D} f(x) ]^2 \dif x}$,
which is equivalent to the norm in~\eqref{eq:sobolev-norm} for $\alpha = 1$.
More generally, the covariance kernel
\begin{equation*}
        K(x, y) = \sum_{k=0}^m \frac{(xy)^k}{(k!)^2} + K_m(x, y) ,
\end{equation*}
where $K_m$ is the $m$ times integrated Brownian motion kernel defined via the recursion
\begin{equation} \label{eq:ibm-kernel}
        K_m(x, y) = \int_0^x \int_0^y K_{m-1}(t, s) \dif t \dif s \quad \text{ and } \quad K_0(x, y) = \min\{x, y\}
\end{equation}
and term $\sum_{k=0}^m (xy)^k / (k!)^2$ serves to remove boundary conditions at the origin, is a Sobolev kernel of order $\nu = m - 1/2$ on $D = [0, 1]$.
See \cite[Sec.\@~3.1]{AdlerTaylor2007},~\cite[Sec.\@~10]{VaartZanten2008}, and~\cite[Sec.\@~1.2]{Wahba1990} for these results.
\Cref{fig:bm-samples} shows samples from integrated Brownian motions.

\begin{remark} \label{rmk:equivalence-classes}
By definition, an RKHS is a space of \emph{functions} while, a priori, the elements of Sobolev spaces are \emph{equivalence classes} of functions that coincide almost everywhere.
The notion of point evaluation, $f(x)$, is meaningless for equivalence classes.
However, the Sobolev embedding theorem~\citep[Thm.\@~4.12]{AdamsFournier2003} ensures that a Sobolev space of order $\alpha$ is continuously embedded in the space of continuous functions if $\alpha > d/2$, in which case each element can be uniquely identified with a continuous function and the space is an RKHS.
The requirement $\alpha > d/2$ is embedded in \Cref{def:sobolev-kernel}.
Note that the role of continuity is but incidental as in general an RKHS can contain discontinuous functions (\cite{SchabackWendland2001} contains an interesting example).
\end{remark}

\section{Asymptotic analysis} \label{sec:results}

How well the observation locations cover the domain $D$ is measured by the \emph{fill-distance}, $h_n$, and the \emph{separation radius}, $q_n$, which are defined as
\begin{equation} \label{eq:fill-distance}
        h_{n} = \sup_{x \in D} \min_{1 \leq i \leq n} \lVert x - x_i \rVert \quad \text{ and } \quad q_n = \frac{1}{2} \min_{1 \leq i \neq j \leq n} \lVert x_i - x_j \rVert .
\end{equation}
The fill-distance is the radius of the largest ball in $D$ that contains none of $x_1, \ldots, x_n$, while the separation radius is half the smallest distance between any two observation locations.
The fill-distance depends on $D$ but the separation radius does not.
When $D$ is convex and has non-empty interior, then $q_n \leq h_{n}$ and any sequence $(x_i)_{i=1}^\infty$ satisfies $h_{n} \gtrsim n^{-1/d}$ and $q_n \lesssim n^{-1/d}$~\cite[Lem.\@~2.1]{PronzatoZhigljavsky2023}.
For example, if $n = (m + 1)^d$ and $x_1, \ldots, x_n$ are the tensor grid formed out of the $m+1$ equispaced points $\{0, 1/m, \ldots, 1\}$ on $[0, 1]$, then
\smash{$h_n = \frac{1}{2}\sqrt{d} m^{-1} = \frac{1}{2}\sqrt{d} (n^{1/d} - 1)^{-1}$} if $D = [0, 1]^d$.
Most of our results assume that the sequence of observation locations is quasi-uniform, which means that $h_n$ and $q_n$ are comparable.

\begin{definition}[Quasi-uniform sequence] \label{def:quasi-uniform}
Let $D \subset \R^d$ be bounded.
A sequence $(x_i)_{i=1}^\infty \subset D$ is \emph{quasi-uniform} if there is $c \geq 1$ such that $c^{-1} q_n \leq h_{n} \leq c \, q_n$ for all $n \geq 1$. 
\end{definition}

A quasi-uniform sequence covers the domain somewhat uniformly. Moreover, $q_n \asymp h_{n} \asymp n^{-1/d}$ if the sequence is quasi-uniform on a bounded open convex $D \subset \R^d$. Note that the empty set is the only open subset of $\R^d$ with empty interior. 
While we consider a sequence of points, as is natural if the observations are obtained sequentially, the definition of quasi-uniformity and our results easily extend to \emph{sequences of point sets} $(\mathcal{X}_n)_{n=1}^\infty$, where $\mathcal{X}_n = (x_{n,k})_{k=1}^n$.

\subsection{Oversmoothing} \label{sec:results-sigma}

Our most general results cover oversmoothing, which refers to $K$ being smoother than $K_0$.
The first result shows that, in expectation, the maximum likelihood and cross-validation scale estimators blow up with a polynomial rate that depends on the extent of oversmoothing by the model.
See \Cref{sec:proofs-results-sigma} for a proof.

\begin{theorem} \label{thm:estimator-asymptotics}
Suppose that $D$ is a bounded open convex set and $K \in \Sob(\nu)$ and $K_0 \in \Sob(\nu_0)$ for $\nu \geq \nu_0 > 0$.
Let $\hat{\sigma}_n^2$ be either the maximum likelihood estimator $\hat{\sigma}_{\ML,n}^2$ in~\eqref{eq:sigma-ML} or the cross-validation estimator \smash{$\hat{\sigma}_{\CV,n}^2$} in~\eqref{eq:sigma-CV}.
If $(x_i)_{i=1}^\infty \subset D$ is quasi-uniform, then
\begin{equation} \label{eq:sigma-E-asymp}
        \E[ \hat{\sigma}_n^2 ] \asymp  n^{2(\nu-\nu_0)/d} 
        \quad \text{ and } \quad 
        \hat{\sigma}_n^2 \asymp_{\mathbb{P}} n^{2(\nu-\nu_0)/d} .
\end{equation}
\end{theorem}

We believe that~\eqref{eq:sigma-E-asymp} holds almost surely but the proof eludes us.
In~\cite[Cor.\@~4.13]{KorteStapff2024} the technique from~\cite[Sec.\@~4.4]{VaartZanten2011} was used to prove an almost sure variant of~\eqref{eq:sigma-E-asymp} for the maximum likelihood estimator when $X_0$ is a Whittle--Matérn process on a closed Riemannian manifold.
As anticipated in \Cref{sec:scale-estimation}, we can use \Cref{thm:estimator-asymptotics} to show that maximum likelihood estimation and cross-validation induce rate-unbiasedness. This is the main result of the article.

\begin{theorem} \label{thm:rate-unbiased}
Suppose that $D$ is a bounded open convex set and $K \in \Sob(\nu)$ and $K_0 \in \Sob(\nu_0)$ for $\nu \geq \nu_0 > 0$.
Let $\hat{\sigma}_n^2$ be either the maximum likelihood estimator $\hat{\sigma}_{\ML,n}^2$ in~\eqref{eq:sigma-ML} or the cross-validation estimator \smash{$\hat{\sigma}_{\CV,n}^2$} in~\eqref{eq:sigma-CV}.
If $(x_i)_{i=1}^\infty \subset D$ is quasi-uniform, then the presumed mean-square error is rate-unbiased for any $p \in (0, \infty)$:
\begin{equation} \label{eq:MSE-unbiased-thm}        
        \E\big[ \, \lVert \MSE_n(\cdot \mid \hat{\sigma}_n) \rVert_{L^p(D)} \, \big] = \E[ \hat{\sigma}_n^2 ] \cdot \lVert \MSE_n \rVert_{L^p(D)} \asymp \lVert \MSE_n^* \rVert_{L^p(D)} .
\end{equation}
\end{theorem}
\begin{proof}
\Cref{thm:MSE-rates-informal,thm:estimator-asymptotics} yield
\begin{equation*}
        \E[ \hat{\sigma}_n^2 ] \cdot  \lVert \MSE_n \rVert_{L^p(D)} \asymp n^{2(\nu-\nu_0)/d} \cdot n^{-2\nu/d} = n^{-2\nu_0/d} \asymp \lVert \MSE_n^* \rVert_{L^p(D)} . \qedhere
\end{equation*}
\end{proof}

\subsection{Undersmoothing} \label{sec:periodic}

The proof of \Cref{thm:estimator-asymptotics} is based on certain escape results for kernel-based interpolation, which state that a kernel interpolant (i.e., the conditional mean $\mu_n$) converges with the optimal rate even if the target function is rougher than expected~\citep{Narcowich2006}.
As such results are not available in sufficient generality if the target function is smoother than expected, we are not able to provide a satisfactory version of \Cref{thm:estimator-asymptotics} that would apply to undersmoothing (i.e., the case $\nu_0 > \nu$).
Something can be said in three special cases, which provide valuable clues to the behaviour of \smash{$\sigma_{\ML,n}^2$} and \smash{$\sigma_{\CV,n}^2$} when the model undersmooths.

First, asymptotics can be derived for the maximum likelihood estimator when there is sufficient undersmoothing.
See \Cref{sec:proofs-results-sigma} for the proof of \Cref{thm:mle-undersmooth-exp}.

\begin{theorem} \label{thm:mle-undersmooth-exp}
Suppose that $D$ is a bounded open convex set and $K \in \Sob(\nu)$ and $H(K_0) \subseteq H_2^\alpha(D)$ for $\nu > 0$ and $\alpha > \nu + d$.
Let $(x_i)_{i=1}^\infty \subset D$ be any sequence.
\begin{enumerate}
        \item[(a)] If $K$ is continuous on $D \times D$, then
        \begin{equation} \label{eq:mle-undersmooth-exp}
        \E[ \sigma_{\ML,n}^2 ] = \frac{1}{n} \tr( \b{K}_n^{-1} \b{K}_{0,n} ) \quad \text{ and } \quad \lim_{n \to \infty} \tr( \b{K}_n^{-1} \b{K}_{0,n} ) \in (0, \infty).
        \end{equation}
        Additionally, $\sigma_{\ML,n}^2 \asymp_\mathbb{P} n^{-1}$.
        \item[(b)] If almost all sample paths of $X_0$ are continuous, then $\sigma_{\ML,n}^2 \asymp n^{-1}$ almost surely.
\end{enumerate}
\end{theorem}

Because the RKHS of $K \in \Sob(\nu)$ is norm-equivalent to a Sobolev space of order $\nu + d/2$, \Cref{thm:mle-undersmooth-exp} states that $\sigma_{\ML,n}^2$ decays in expectation as $n^{-1}$ whenever the truth is smoother than the model by at least order $d/2$.
The covariance $K_0$ satisfies the assumption in \Cref{thm:mle-undersmooth-exp} if it is a Sobolev kernel of order $\nu_0 > \nu + d/2$.
Note that the interval $\nu_0 \in (\nu, \nu + d/2]$ is not covered by either of Theorems~\ref{thm:estimator-asymptotics} or~\ref{thm:mle-undersmooth-exp}.
In addition to Sobolev kernels, the theorem applies to a variety of commonly used non-Sobolev kernels, such as the Gaussian covariance in~\eqref{eq:gaussian-kernel}.
Since the spectral density of the Gaussian covariance is a Gaussian function and thus decays faster than any polynomial, its RKHS is contained in every Sobolev space by~\eqref{eq:RKHS-fourier} and~\eqref{eq:Bessel-norm}.
The trace limit in~\eqref{eq:mle-undersmooth-exp} is in fact the trace of a certain operator between the RKHSs of $K$ and $K_0$~\citep[Prop.\@~4.5]{LukicBeder2001}.
The properties of this operator control whether or not the samples of $X_0$ are contained in $H(K)$; see~\citep{Driscoll1973} and \citep{LukicBeder2001}, as well as the proof of \Cref{thm:mle-undersmooth-exp} in \Cref{sec:proofs-results-sigma}.

\begin{figure}
        \centering
        \includegraphics[width=\textwidth]{\figpath 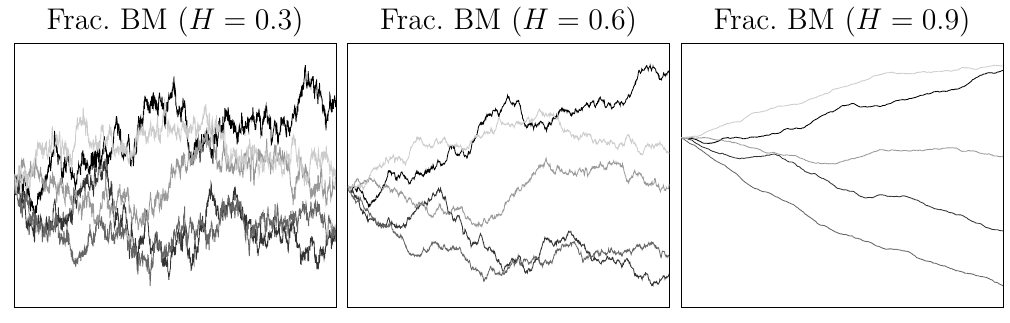}
        \caption{Samples from the fractional Brownian motion with index $H \in \{0.3, 0.6, 0.9\}$. Covariance kernels of these Gaussian processes are given in~\eqref{eq:FBM-kernel}.}
        \label{fig:fbm-samples}
\end{figure}

Second, Naslidnyk et al.\@~\citep[Thms.\@~11 and~12]{Naslidnyk2024} have proved an undersmoothing theorem under the assumption that $K(x, y) = \min\{x, y\}$ is the Brownian motion kernel on $D = [0, 1]$ and $X_0$ is either the fractional Brownian motion with the Hurst index $H \in (0, 1)$ or its integral.
The covariances of the fractional Brownian motion and its integral are
\begin{equation} \label{eq:FBM-kernel}
        K_{0,H}(x, y) = \frac{1}{2} ( \lvert x \rvert^{2H} + \lvert y \rvert^{2H} - \lvert x-y\rvert^{2H} )
\end{equation}
and 
\begin{equation} \label{eq:iFBM-kernel}
  \begin{split}
  K_{1,H}(x, y) &= \int_0^x \int_0^{y} K_{0,H}(t, t') \dif t \dif t' \\
  &= \frac{1}{2H'} \bigg( y x^{H'} + x y^{H'} - \frac{1}{H' + 1} \big( x^{H'+1} + y^{H'+1} - \lvert x - y \rvert^{H'+1} \big) \bigg) ,
  \end{split}
\end{equation}
where $H' = 2H+1$.
Recall that the index $H = 1/2$ recovers the Brownian motion, in which case the model is well-specified.
\Cref{fig:fbm-samples} shows samples from fractional Brownian motions.

\begin{theorem}[Thms.\@~11 and~12 in \citep{Naslidnyk2024}] \label{thm:naslidnyk}
Suppose that $D = [0, 1]$, that $(x_i)_{i=1}^\infty \subset [0, 1]$ is quasi-uniform, and that $X$ is the Brownian motion.
If $X_0$ is the fractional Brownian motion with Hurst index $H \in (0, 1)$ and covariance~\eqref{eq:FBM-kernel}, then
\begin{equation} \label{eq:BM-FBM-sigma}
        \E[ \hat{\sigma}_{\ML,n}^2 ] \asymp \E[ \hat{\sigma}_{\CV,n}^2] \asymp n^{1-2H} . 
\end{equation}
If $X_0$ is the integrated fractional Brownian motion with $H \in (0, 1)$ and covariance~\eqref{eq:iFBM-kernel}, then
\begin{equation*}
        \E[ \hat{\sigma}_{\ML,n}^2 ] \asymp n^{-1} \quad \text{ and } \quad \E[ \hat{\sigma}_{\CV,n}^2] \asymp 
        \begin{cases}
                n^{-1-2H} &\text{ if } \quad H < 1/2, \\
                n^{-2} &\text{ if } \quad H \geq 1/2 .
        \end{cases}
\end{equation*}
All rates hold also in the sense of~\eqref{eq:asymp-in-prob}.
\end{theorem}

Note that the rate in~\eqref{eq:BM-FBM-sigma} is constant if $H = 1/2$.
\Cref{thm:naslidnyk} shows that scale estimation may be effective in combatting not only oversmoothing but also undersmoothing. 
Moreover, the theorem suggests that cross-validation should be preferred over maximum likelihood estimation because less undersmoothing is required to saturate the latter (a rigorous justification is provided by \Cref{thm:rate-unbiased-periodic}).
That is, \smash{$\hat{\sigma}_{\ML,n}^2$} decays as $n^{-1}$ for any integrated fractional Brownian motion while the rate for \smash{$\hat{\sigma}_{\CV,n}^2$} depends on the Hurst index up to $H = 1/2$.
Note that \Cref{thm:mle-undersmooth-exp} contains a general saturation result for \smash{$\hat{\sigma}_{\ML,n}^2$} that applies to almost any kernel.
\Cref{thm:naslidnyk} does not generalise easily because its proof uses explicit expressions for the scale estimator unique to the Brownian motion prior.

\begin{figure}
        \centering
        \includegraphics[width=\textwidth]{\figpath 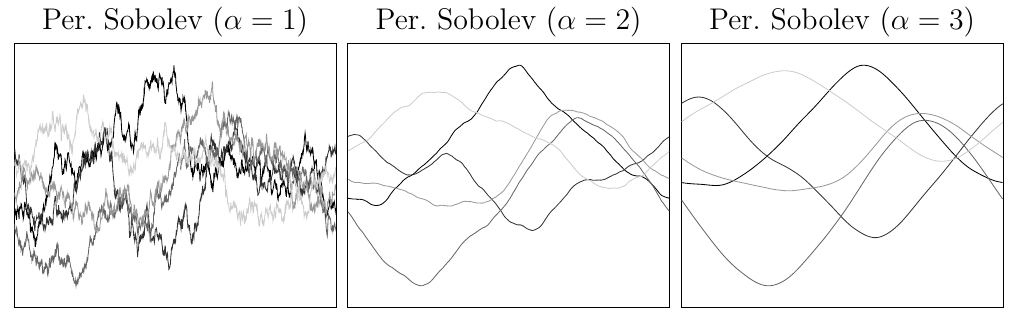}
        \caption{Samples from Gaussian processes with the periodic covariance kernel in~\eqref{eq:korobov-kernel} with $\alpha \in \{1, 2, 3\}$.}
        \label{fig:periodic-samples}
\end{figure}

Third, we are able to prove a limited generalisation of \Cref{thm:naslidnyk} for certain periodic priors and processes.
We refer to Appendix~A.1 in \citep{NovakWozniakowski2008} and Section~2.1 in~\citep{DickKritzer2022} for the following results.
Let $\mathrm{i}$ be the imaginary unit and $\varphi_k(x) = e^{2\pi \mathrm{i} k x}$.
On $D = [0, 1]$, the periodic Sobolev kernel (or Korobov kernel) of order $\alpha > 1/2$ is given by the series expansion
\begin{equation} \label{eq:korobov-kernel-alpha}
        \begin{split}
        K(x, y) = 1 + \sum_{k \neq 0} \lvert k \rvert^{-2 \alpha} e^{2\pi \mathrm{i} k x} \overline{e^{2\pi \mathrm{i} k y}} &= 1 + \sum_{k \neq 0} \lvert k \rvert^{-2 \alpha} e^{2\pi \mathrm{i} k (x-y)} \\
        &= 1 + 2 \sum_{k=1}^\infty k^{-2 \alpha} \cos(2\pi k(x-y)) .
        \end{split}
\end{equation}
For $\alpha \in \N$, the kernel can be written as
\begin{equation} \label{eq:korobov-kernel}
        \begin{split}
        K(x, y) = 1 + (-1)^{\alpha+1} (2\pi)^{2\alpha} \frac{\mathrm{B}_{2\alpha}( \lvert x - y \rvert )}{(2\alpha)!} ,
        \end{split}
\end{equation}
where $\mathrm{B}_{2\alpha}$ is the Bernoulli polynomial of degree $2\alpha$.
The RKHS of $K$ is the periodic Sobolev space of order $\alpha$.
If $\alpha \in \N$, the RKHS consists of those functions in the Sobolev space $H^\alpha([0,1])$ whose derivatives up to order $\alpha-1$ are periodic and its norm is equivalent to the Sobolev norm.
Works on periodic Sobolev kernels in the statistics literature include~\cite{Chen2021, Petit2023, SanzAlonso2024}.
\Cref{fig:periodic-samples} shows samples from periodic Sobolev processes.
The following theorem is a generalisation of \Cref{thm:naslidnyk} to periodic Sobolev kernels.
See \Cref{sec:proofs-results-sigma} for a proof.

\begin{theorem} \label{thm:periodic}
Suppose that $D = [0, 1]$ and that $K$ and $K_0$ are periodic Sobolev kernels of orders $\alpha \in \N$ and $\alpha_0 \in \R$ such that $\alpha_0 \geq \alpha > 1/2$.
If $(x_i)_{i=1}^\infty \subset D$ is quasi-uniform, then
\begin{equation*}
        \E[ \hat{\sigma}_{\ML,n}^2 ] \asymp
        \begin{cases}
                n^{2(\alpha - \alpha_0)} &\text{ if } \quad \alpha_0 < \alpha + 1/2 , \\
                n^{-1} \log n &\text{ if } \quad \alpha_0 = \alpha + 1/2, \\
                n^{-1} &\text{ if } \quad \alpha_0 > \alpha + 1/2 
        \end{cases}
\end{equation*}
and
\begin{equation*}
        \E[ \hat{\sigma}_{\CV,n}^2 ] \asymp 
        n^{2(\alpha - \alpha_0)} \quad \text{ if } \quad \alpha_0 \leq 2\alpha.
\end{equation*}
All rates hold also in the sense of~\eqref{eq:asymp-in-prob}.
\end{theorem}

We believe that the condition $\alpha \in \N$ can be removed but have been unable to locate an appropriate result on approximation in periodic Sobolev spaces that would allow this.

\begin{remark}
The curious polylogarithmic rate $n^{-1} \log n$ for the maximum likelihood estimator in \Cref{thm:periodic} is nowhere to be seen in \Cref{thm:naslidnyk}.
The case $\alpha_0 = \alpha + 1/2$ would correspond to $H = 1$ in \Cref{thm:naslidnyk}.
However, this case is not covered by \Cref{thm:naslidnyk} because $H = 1$ does not correspond to a valid fractional Brownian motion.
\end{remark}

The following rate-unbiasedness theorem shows that the range of regularities for which cross-validation yields rate-unbiasedness is significantly larger than for maximum likelihood when a periodic model undersmooths.
Cross-validation ensures rate-unbiasedness for all $\alpha_0 \in [\alpha, 2\alpha]$ while maximum likelihood is limited to $\alpha_0 \in [\alpha, \alpha + 1/2]$.
In this sense one should prefer cross-validation over maximum likelihood if the model is misspecified.

\begin{theorem} \label{thm:rate-unbiased-periodic}
Suppose that $D = [0, 1]$ and that $K$ and $K_0$ are periodic Sobolev kernels of orders $\alpha \in \N$ and $\alpha_0 \in \R$ such that $\alpha_0 \geq \alpha > 1/2$.
If $(x_i)_{i=1}^\infty \subset D$ is quasi-uniform and $p \in (0, \infty)$, then
\begin{equation*}
        \E\big[ \, \lVert \MSE_n(\cdot \mid \hat{\sigma}_{\ML,n}) \rVert_{L^p(D)} \, \big] \asymp \lVert \MSE_n^* \rVert_{L^p(D)} \quad \text{ if and only if } \quad \alpha_0 < \alpha + 1/2
\end{equation*}        
and
\begin{equation*}
        \E\big[ \, \lVert \MSE_n(\cdot \mid \hat{\sigma}_{\CV,n}) \rVert_{L^p(D)} \, \big] \asymp \lVert \MSE_n^* \rVert_{L^p(D)} \quad \text{ if } \quad \alpha_0 \leq 2\alpha .
\end{equation*}
\end{theorem}
\begin{proof}
The claims follow from \Cref{thm:periodic} after applying \Cref{thm:MSE-bounds-periodic}, which yields
\begin{equation*}
        \E[ \hat{\sigma}_n^2 ] \cdot  \lVert \MSE_n \rVert_{L^p(D)} \asymp \E[ \hat{\sigma}_n^2 ] \cdot n^{-2\alpha + 1} \quad \text{ and } \quad \lVert \MSE_n^* \rVert_{L^p(D)} \asymp n^{-2\alpha_0 + 1}. \qedhere
\end{equation*}
\end{proof}

\section{Numerical examples and smoothness estimation} \label{sec:numerical}

This section contains two numerical examples that validate the theory of \Cref{sec:results} and demonstrate that the smoothness of $X_0$ can be estimated from the rate of growth of a scale estimator via \Cref{thm:estimator-asymptotics}.

\subsection{Smoothness estimation} \label{sec:smoothness-estimation}

\Cref{thm:estimator-asymptotics} states that the maximum likelihood and cross-validation estimators blow up as \smash{$n^{2(\nu-\nu_0)/d}$} if $K \in \Sob(\nu)$ and \smash{$K_0 \in \Sob(\nu_0)$} for $\nu \geq \nu_0 > 0$.
This suggests that the smoothness $\nu_0$ of $X_0$ can be estimated by computing a scale estimator for a number of different $n$ and fitting a polynomial to the resulting data.
The leading order of the fitted polynomial determines the smoothness.
The following approach appears simplest:
\begin{enumerate}
\item Select $\nu > 0$. This parameter should satisfy $\nu \geq \nu_0$.
\item Select $0 \leq n_1 < \cdots < n_m \leq n$ and compute $\hat{\sigma}_{n_i}^2 \in \{\hat{\sigma}_{\ML,n_i}^2, \hat{\sigma}_{\CV,n_i}^2\}$ for $i \in \{1, \ldots m\}$.
\item Fit a linear polynomial to the resulting data on logarithmic scale. That is, solve the linear least-squares problem
\begin{equation*}        
        (\hat{a}_n, \hat{\beta}_n) = \argmin_{a, \beta \in \R} \, \sum_{i=1}^m ( a + \beta \log n_i - \log \hat{\sigma}_{n_i}^2 )^2 .
\end{equation*}
On the linear scale, this gives the fit $e^{\hat{a}_n} n^{\hat{\beta}_n}$ to the data $\{(n_i, \hat{\sigma}_{n_i}^2)\}_{i=1}^m$.
\item Because the scale estimators are expected to behave as $n^{2(\nu-\nu_0)/d}$ if $\nu \geq \nu_0$, take
\begin{equation} \label{eq:nu-estimator}
        \hat{\nu}_n = \nu - d \hat{\beta}_n / 2 
\end{equation}
as a smoothness estimate.
\end{enumerate}

Countless minor and obvious variations are possible.
For example, fitting a polynomial of the form $a_0 + a_1 n^\beta$ might yield better estimates when $n$ is small at the cost of making the optimisation problem non-linear.
A major disadvantage of the method is that it requires oversmoothing ($\nu \geq \nu_0$), which is difficult to guarantee in practice.
The following proposition shows that the method recovers the true smoothness.
The proof is given in \Cref{sec:proofs-numerical}

\begin{proposition} \label{prop:smoothness:LS}
Consider the setting of \Cref{thm:estimator-asymptotics} and the smoothness estimator $\hat{\nu}_n$ in~\eqref{eq:nu-estimator} for a fixed $m \geq 2$.
Assume that $n_1 = n_1(n) \to \infty, \ldots , n_m = n_m(n) \to \infty$ and $n_m(n) / n_1(n) \to \infty$ as $n \to \infty$. Then $\hat{\nu}_n \to \nu_0$ in probability.
\end{proposition}

\subsection{Ruzsa's sequence}

Though hidden by the asymptotic notation, the bounds in \Cref{thm:estimator-asymptotics} depend on the ratio $h_n / q_n$ of the fill-distance and separation radius defined in~\eqref{eq:fill-distance}.
The dependency is explicit in \Cref{thm:arcangeli} that we use to prove \Cref{thm:estimator-asymptotics}.
To make it easier to validate \Cref{thm:estimator-asymptotics} we thus want a sequence of points for which $h_n$ and $q_n$ vary as smoothly as possible.
If we elected to work with a non-nested sequence of point sets, by which we mean sets of points $\mathcal{X}_n = (x_{n,k})_{k=1}^n$ such that $\mathcal{X}_n \not\subset \mathcal{X}_{n+1}$, we could simply take sets of equispaced points.
However, in practice one either obtains data sequentially or, when using the smoothness estimation method from~\Cref{sec:smoothness-estimation}, starts from a given point set and constructs a sequence of nested subsets.

The well-known van der Corput sequence $(0, 1, \tfrac{1}{2}, \tfrac{1}{4}, \tfrac{3}{4}, \tfrac{1}{8}, \ldots )$ in base $2$ on $D = [0, 1]$ is an obvious candidate for a uniform sequence~\cite[Def.\@~3.2]{Niederreiter1992}.
However, $h_n$ and $q_n$ of the van der Corput sequence behave non-smoothly.
For each $k \geq 0$, we have  $h_n = 2^{-k-1}$ for $n \in \{ 2^{k}, \ldots, 2^{k+1} - 1\}$ and $q_n = 2^{-k-2}$ for $n \in \{ 2^{k} + 1, \ldots, 2^{k+1}\}$
One can do much better by using \emph{Ruzsa's sequence}~\cite[p.\@~154]{Niederreiter1992} given by
\begin{equation*} 
        x_1 = 1 \quad \text{ and } \quad x_k = \{ \log_2(2n - 3) \} \quad \text{ for } \quad k \geq 2,
\end{equation*}
where $\{x\} = x - \lfloor x \rfloor$ is the fractional part.
It is straightforward to compute that 
\begin{equation*}
        h_n = \frac{1}{\log 4} n^{-1} + O(n^{-2}), \quad q_n = \frac{1}{2 \log 4} n^{-1} + O(n^{-2}), \quad \text{ and } \quad \frac{h_n}{q_n} = 2 + O(n^{-1}) 
\end{equation*}
for Rusza's sequence.
Ruzsa's sequence is known to have the smallest possible covering constant, $\limsup_{n \to \infty} n \cdot h_n = 1/\log 4$, among nested sequences~\cite[Thm.\@~6.7]{Niederreiter1992}.

\begin{figure}[t]
        \centering
        \includegraphics[width=\textwidth]{\figpath 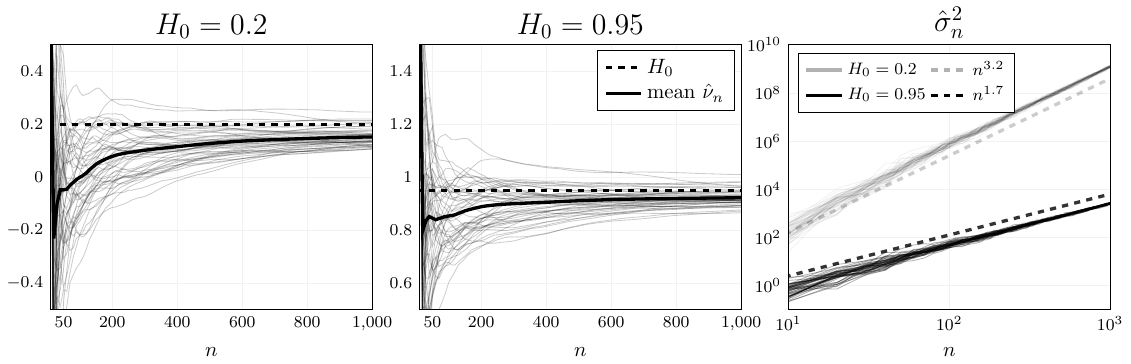}
        \caption{Smoothness estimator $\hat{\nu}_n$ from \Cref{sec:smoothness-estimation} (left and middle) and the maximum likelihood scale estimator $\hat{\sigma}_n^2 = \hat{\sigma}_{\ML,n}^2$ (right) for 50 samples from $X_0$ in the univariate Case~1 described in \Cref{sec:numerical-results}.}
        \label{fig:num-ex-1}
\end{figure}

\begin{figure}[t]
        \centering
        \includegraphics[width=\textwidth]{\figpath 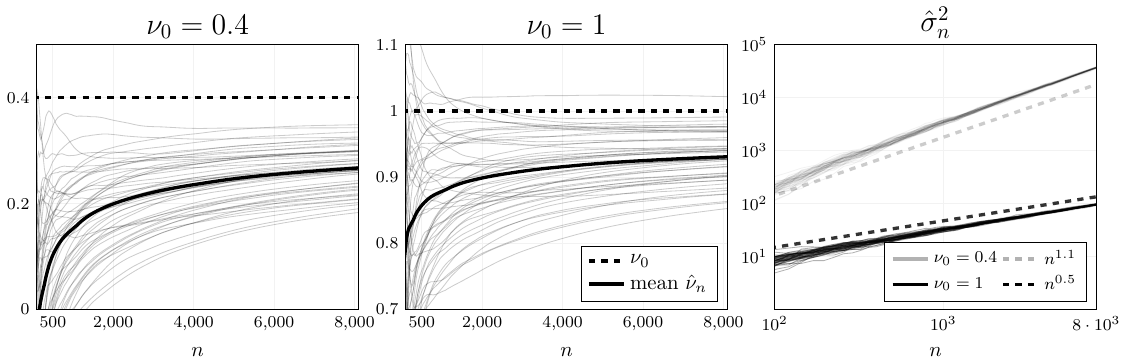}
        \caption{Smoothness estimator $\hat{\nu}_n$ from \Cref{sec:smoothness-estimation} (left and middle) and the maximum likelihood scale estimator $\hat{\sigma}_n^2 = \hat{\sigma}_{\ML,n}^2$ (right) for 50 samples from $X_0$ in the bivariate Case~2 described in \Cref{sec:numerical-results}. Note that, due to the curse of dimensionality, \smash{$n = 8,\!100 = 90^2$} roughly corresponds to $n = 90$ in \Cref{fig:num-ex-1}.}
        \label{fig:num-ex-2}
\end{figure}

\subsection{Setting and results for Ruzsa's sequence} \label{sec:numerical-results}

We study the behaviour of the maximum likelihood estimator $\hat{\sigma}_n^2 = \hat{\sigma}_{\ML,n}^2$ and the corresponding smoothness estimator $\hat{\nu}_n$ in~\eqref{eq:nu-estimator}:

\begin{enumerate}
        \item We take $d = 1$ and $D = [0, 1]$. 
        The prior covariance is a Matérn with smoothness $\nu = 1.8$ and correlation length $\lambda = 1$, while the true process $X_0$ is a fractional Brownian motion with Hurst index (a) $H_0 = 0.2$ or (b) $H_0 = 0.95$.        
        We use Ruzsa's sequence and compute $\hat{\sigma}_n^2$ for every 10th $n$ up to $n = 1,\!000$.
        Based on these scale estimators up to a given $n$ we compute $\hat{\nu}_n$.
        Although fractional Brownian motion kernels do not appear to be Sobolev kernels (see~\cite{BartonPoor1988} for an RKHS characterisation) and \Cref{thm:estimator-asymptotics} is not strictly speaking applicable, we nevertheless make the identification $H_0 = \nu_0$. Recall that the covariance of the Brownian motion (i.e., $H_0 = 0.5$) is a Sobolev kernel of order $\nu_0 = 0.5$ up to a boundary condition at the origin.
        We thus expect to observe the rates
        \begin{equation} \label{eq:rates-ex1}
                \E[ \hat{\sigma}_n^2 ] \asymp  n^{2(\nu-H_0)/d} = n^{3.2} \quad \text{ and } \quad \E[ \hat{\sigma}_n^2 ] \asymp  n^{2(\nu-H_0)/d} = n^{1.7}
        \end{equation}
        for $H_0 = 0.2$ and $H_0 = 0.95$, respectively.

        \item We take $d = 2$ and $D = [0, 1]^2$. 
        The prior covariance is a Matérn with smoothness $\nu = 1.5$ and correlation length $\lambda = 1.5$, while the true process $X_0$ is a Matérn with $\lambda = 1$ and (a) $\nu_0 = 0.4$ or (b) $\nu_0 = 1$.        
        We use product grids formed using the $\tilde{n} \in \{3, 4, \ldots, 90\}$ first points of Ruzsa's sequence and compute $\hat{\sigma}_n^2$ for $n = \tilde{n}^2 \in \{3^2, 4^2, \ldots, 90^2 = 8,\!100\}$.
        Based on these scale estimators up to a given $n$ we compute $\hat{\nu}_n$.
        In this case all assumptions of \Cref{thm:estimator-asymptotics} are satisfied and we expect to observe the rates
        \begin{equation} \label{eq:rates-ex2}
                \E[ \hat{\sigma}_n^2 ] \asymp  n^{2(\nu-\nu_0)/d} = n^{1.1} \quad \text{ and } \quad \E[ \hat{\sigma}_n^2 ] \asymp  n^{2(\nu-\nu_0)/d} = n^{0.5}
        \end{equation}
        for $\nu_0 = 0.4$ and $\nu_0 = 1$, respectively.        
\end{enumerate}

The results are shown in \Cref{fig:num-ex-1,fig:num-ex-2} for 50 samples from $X_0$.
The rates in~\eqref{eq:rates-ex1} and~\eqref{eq:rates-ex2} predicted by \Cref{thm:estimator-asymptotics} are clearly observed and the estimator $\hat{\nu}_n$ converges.
The convergence is quite slow and the estimator appears to be biased.
The convergence of the standard maximum likelihood smoothness estimator \smash{$\hat{\nu}_{\ML, n} = \argmin_{\nu > 0} \{ \b{y}_n^\T \b{K}_{\nu, n}^{-1} \b{y}_n + \log \det \b{K}_{\nu,n} \}$}, which is not shown here, is much faster in these examples.
However, to compute this estimator one has to repeatedly form and perform linear algebra with the covariance matrix $\b{K}_{\nu,n}$ for different $\nu$, which is computationally expensive.
The scale estimation based smoothness estimator $\hat{\nu}_n$ may be a viable option whenever cheap and rough smoothness estimates suffice.

\subsection{Results for other sequences}

\begin{figure}
        \centering
        \includegraphics[width=0.7\textwidth]{\figpath 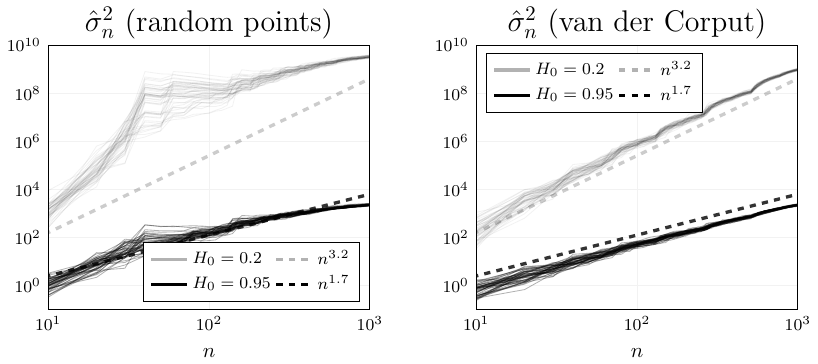}
        \caption{The maximum likelihood scale estimator $\hat{\sigma}_n^2 = \hat{\sigma}_{\ML,n}^2$ for 50 samples from $X_0$ in the univariate Case~1 described in \Cref{sec:numerical-results} with the exception that here the points are drawn randomly from the uniform distribution on $D = [0, 1]$ (left) and from the van der Corput sequence (right). The estimators behave more erratically than for Rusza's sequence (\Cref{fig:num-ex-1,fig:num-ex-2}).}
        \label{fig:num-ex-other}
\end{figure}

For sequences other than Rusza's the results are not quite as nice.
\Cref{fig:num-ex-other} shows the behaviour of $\hat{\sigma}_n^2$ for 50 samples in Case~1, where $X$ is a Matérn, $X_0$ a fractional Brownian motion, and $D = [0, 1]$.
But now we use either (a) the first $n$ points from a sequence of 1,000 points drawn randomly from the uniform distribution on $[0, 1]$ or (b) the van der Corput sequence.
The scale estimators behave much more erratically than for Rusza's sequence, particularly when the points are random.
The oscillations for the van der Corput sequence correspond to powers of two: the sequence is equispaced for $n = 2^k$.

\section{Conclusion}

We have proved that computationally tractable scale estimation is a powerful tool to combat smoothness misspecification in Gaussian process modelling.
If the model oversmooths, the true mean-square error and that presumed by the model will decay with the same rate in expectation and probability if the scale parameter is set with maximum likelihood or cross-validation.
If the model undersmooths (but not too much), cross-validation outperforms maximum likelihood estimation in this sense.
There are two significant missing pieces.
First, our results hold in expectation. 
We believe that the results are valid almost surely but do not presently know how to prove this.
Second, our undersmoothing results concern periodic Sobolev kernels rather than popular Matérns.

\section{Proofs}

All proofs longer than a few lines are collected in this section.

\subsection{Kernel interpolation} \label{sec:kernel-interpolation}

The proofs are based on error estimates for kernel interpolants in Sobolev spaces.
Throughout this section $K \colon D \times D \to \R$ is a positive-definite kernel.
Let $x_1, \ldots, x_n \in D$ be pairwise distinct points and $f \colon D \to \R$ a function.
The \emph{kernel interpolant} to $f$ at these points is the function $I_n f$ given by
\begin{equation} \label{eq:kernel-interpolant}
        (I_n f)(x) = \b{k}_n(x)^\T \b{K}_n^{-1} \b{f}_n,
\end{equation}
where $\b{k}_n(x) \in \R^n$ and $\b{K}_n \in \R^{n \times n}$ are as in~\eqref{eq:GP-mean-cov} and $\b{f}_n = (f(x_1), \ldots, f(x_n)) \in \R^n$ collects evaluations of $f$.
The kernel interpolation operator $I_n$ is the mapping $f \mapsto \b{k}_n(\cdot)^\T \b{K}_n^{-1} \b{f}_n$.
We see that the kernel interpolant is nothing but the posterior mean $\mu_n$ with the observations $\b{y}_n = \b{f}_n$.
The kernel interpolant is constructed as the unique function in the linear span of the kernel translates $K(x_1, \cdot), \ldots, K(x_n, \cdot)$ that interpolates $f$ at the points $x_1, \ldots, x_n$.
That is, $(I_n f)(x) = \sum_{i=1}^n a_i K(x, x_i)$ for coefficients $a_i$ that solve the linear system of equations
\begin{equation*}
\begin{bmatrix}
K(x_1, x_1) & \cdots & K(x_1, x_n) \\
\vdots & \ddots & \vdots \\
K(x_n, x_1) & \cdots & K(x_n, x_n) 
\end{bmatrix}
\begin{bmatrix}
        a_1 \\ \vdots \\ a_n
\end{bmatrix}
=
\begin{bmatrix}
        f(x_1) \\ \vdots \\ f(x_n)
\end{bmatrix} .
\end{equation*}
Solving this equations yields~\eqref{eq:kernel-interpolant}.

Recall the notion of an RKHS from \Cref{sec:sobolev}.
Let $R \colon D \times D \to \R$ be another positive-definite kernel with an RKHS $H(R)$.
Let $x \in D$. 
The \emph{worst-case error} in $H(R)$ of the kernel interpolant is defined as the largest interpolation error for functions in the unit ball of $H(R)$:
\begin{equation} \label{eq:wce-def}
        e_n(x; R) = \sup_{ \lVert f \rVert_{H(R)} \leq 1} \lvert f(x) - (I_n f)(x) \rvert . 
\end{equation}
Note that the interpolant is constructed with the kernel $K$, which need not equal $R$.
The worst-case error has a well-known closed-form expression that can be derived from, for example, Section~10.2 in \citep{NovakWozniakowski2010}.
The proof is not difficult, so we provide it for completeness.

\begin{proposition}
Let $I_n$ be the kernel interpolation operator given by~\eqref{eq:kernel-interpolant}.
Then
\begin{equation} \label{eq:wce-explicit}
        e_n(x; R) = \sqrt{\smash[b]{ R(x, x) - 2 \b{r}_n(x) \b{K}_n^{-1} \b{k}_n(x) + \b{k}_n(x)^\T \b{K}_n^{-1} \b{R}_n \b{K}_n^{-1} \b{k}_n(x) }} ,
\end{equation}
where $\b{r}_n(x) \in \R^n$ and $\b{R}_n \in \R^{n \times n}$ are defined analogously to $\b{k}_n(x)$ and $\b{K}_n$.
\end{proposition}
\begin{proof}
Fix $x \in D$ and write the kernel interpolant as $(I_n f)(x) = \sum_{i=1}^n c_i f(x_i) = \b{c}^\T \b{f}_n$, where $\b{c} = (c_1, \ldots, c_n) = \b{K}_n^{-1} \b{k}_n(x)$.
Using the reproducing property in~\eqref{eq:reproducing-property} we write
\begin{equation*}
        e_n(x; R) = \sup_{ \lVert f \rVert_{H(R)} \leq 1} \bigg\lvert f(x) - \sum_{i=1}^n c_i f(x_i) \bigg\rvert = \sup_{ \lVert f \rVert_{H(R)} \leq 1} \bigg\lvert \bigg\langle f , R(\cdot, x)  - \sum_{i=1}^n c_i R(\cdot, x_i) \bigg\rangle_{\!\!H(R)} \bigg\rvert .
\end{equation*}
The Cauchy--Schwarz inequality yields $e_n(x; R) \leq \lVert R(\cdot, x)  - \sum_{i=1}^n c_i R(\cdot, x_i) \rVert_{H(R)}$.
By the reproducing property and the symmetry of $R$, the square of the norm on the right-hand side is
\begin{equation*}
\begin{split}        
        \langle R(\cdot, x), R(\cdot, x) \rangle_{H(R)} - 2 \sum_{i=1}^n \langle R(\cdot, x) &,  c_i R(\cdot, x_i) \rangle_{H(R)} + \sum_{i,j=1}^n \langle c_i R(\cdot, x_i) , c_j R(\cdot, x_j) \rangle_{H(R)} \\
        ={}& R(x, x) - 2 \sum_{i=1}^n c_i R(x, x_i) + \sum_{i,j=1}^n c_i c_j R(x_i, x_j) \\
        ={}& R(x, x) - 2 \b{r}_n(x)^\T \b{c} + \b{c}^\T \b{R}_n \b{c} .
\end{split}
\end{equation*}
Therefore $e_n(x; R) \leq ( R(x, x) - 2 \b{r}_n(x)^\T \b{c} + \b{c}^\T \b{R}_n \b{c} )^{1/2}$.
To see that this is an equality, select the function
\begin{equation*}
        f = \frac{R(\cdot, x) - \sum_{i=1}^n c_i R(\cdot, x_i)}{\lVert R(\cdot, x) - \sum_{i=1}^n c_i R(\cdot, x_i) \rVert_{H(R)} },
\end{equation*}
which has unit norm in $H(R)$, and compute
\begin{equation*}
        \begin{split}
        e_n(x; R) \geq \bigg\lvert f(x) - \sum_{i=1}^n c_i f(x_i) \bigg\rvert &= \frac{R(x, x) - 2 \b{r}_n(x)^\T \b{c} + \b{c}^\T \b{R}_n \b{c}}{\lVert R(\cdot, x) - \sum_{i=1}^n c_i R(\cdot, x_i) \rVert_{H(R)}} \\
        &= \sqrt{ \smash[b]{R(x, x) - 2 \b{r}_n(x)^\T \b{c} + \b{c}^\T \b{R}_n \b{c} }}.
        \end{split}
\end{equation*}
Inserting  $\b{c} = \b{K}_n^{-1} \b{k}_n(x)$ yields the claim.
\end{proof}

Two things should now be observed.
By selecting $R = K_0$ in~\eqref{eq:wce-explicit} we see that the squared worst-case error equals the true mean-square error in~\eqref{eq:MSE-true}:
\begin{equation} \label{eq:wce-mse-1}
        e_n(x; K_0)^2 = \MSE_n^*(x) .
\end{equation}
If we set $R = K$, the expression for the worst-case error simplifies and we obtain
\begin{equation} \label{eq:wce-mse-2}
        e_n(x; K)^2 = K(x, x) - \b{k}_n(x)^\T \b{K}_n^{-1} \b{k}_n(x) = V_n(x) = \MSE_n(x),
\end{equation}
where $V_n(x)$ and $\MSE_n(x)$ are the posterior variance and mean-square error presumed by the model in~\eqref{eq:posterior-var} and~\eqref{eq:MSE}.

\subsection{Bounds on mean-square errors}

Let $p \in [1, \infty]$. We use the notation
\begin{equation} \label{eq:wce-p}
        e_{n, p}(R) = \sup_{ \lVert f \rVert_{H(R)} \leq 1} \lVert f - I_n f \rVert_{L^p(D)}
\end{equation}
for worst-case error measured in the $L^p(D)$-norm.
The following theorem is a consequence of Corollary~4.1 in \citep{Arcangeli2007} and the results in \citep{Narcowich2006}.
In the context of Gaussian processes this theorem is Theorem~1 of \citep{Wynne2021} with $q = p$, $s = 0$, $\tau_f = \tau + d/2$, and $\tau_k^- = \tau_k^+ = \nu + d/2$.

\begin{theorem} \label{thm:arcangeli}
Suppose that $D$ is a bounded open convex set and $K \in \Sob(\nu)$ and $R \in \Sob(\tau)$ for $\nu \geq \tau > 0$.
Let $p \in [1, \infty]$ and $(x)_+ = \max\{0, x\}$.
Then
\begin{equation*}
        e_{n, p}(R) \lesssim \bigg( \frac{h_{n}}{q_n} \bigg)^{\nu-\tau} h_{n}^{\tau + d/2 - d(1/2 - 1/p)_+}.
\end{equation*}
If the sequence $(x_i)_{i=1}^\infty \subset D$ is quasi-uniform, then
\begin{equation*}
        e_{n, p}(R) \lesssim n^{-\tau/d + 1/2 - (1/2 - 1/p)_+} .
\end{equation*}
\end{theorem}

\Cref{thm:arcangeli} and the equivalences above yield rates for mean-square errors.

\begin{corollary} \label{thm:mse-upper-bound}
Suppose that $D$ is a bounded open convex set and $K \in \Sob(\nu)$ and $K_0 \in \Sob(\nu_0)$ for $\nu \geq \nu_0 > 0$.
Then
\begin{equation*}
        \sup_{x \in D} \MSE_n(x) \lesssim h_{n}^{2\nu} \quad \text{ and } \quad \sup_{x \in D} \MSE_n^*(x) \lesssim \bigg( \frac{h_{n}}{q_n} \bigg)^{2(\nu-\nu_0)} h_{n}^{2\nu_0}.
\end{equation*}
If the sequence $(x_i)_{i=1}^\infty \subset D$ is quasi-uniform, then
\begin{equation*}
        \sup_{x \in D} \MSE_n(x) \lesssim n^{-2\nu/d} \quad \text{ and } \quad \sup_{x \in D} \MSE_n^*(x) \lesssim n^{-2\nu_0/d} .
\end{equation*}
\end{corollary}
\begin{proof}
Observe that
\begin{equation*}
        \sup_{x \in D} \MSE_n^*(x) = \sup_{x \in D} e_n(x; K_0)^2 \leq e_{n,\infty}(K_0)^2 
\end{equation*}
and
\begin{equation*}
        \sup_{x \in D} \MSE_n(x) = \sup_{x \in D} e_n(x; K)^2 \leq e_{n,\infty}(K)^2 
\end{equation*}
by~\eqref{eq:wce-def} and~\eqref{eq:wce-mse-1}--\eqref{eq:wce-p}.
The claims thus follow from \Cref{thm:arcangeli} with $p = \infty$ and $R = K_0$.
\end{proof}

In the following theorem we adopt slightly more general notation and use $\MSE_{\mathcal{X}}$ and $\MSE_{\mathcal{X}}^*$ to denote mean-square errors when $\mathcal{X} \subseteq D$ is a finite set of pairwise distinct observation locations.
That is, $\MSE_{\mathcal{X}} = \MSE_n$ and $\MSE_{\mathcal{X}}^* = \MSE_n^*$ if $\mathcal{X} = \{x_i\}_{i=1}^n$.

\begin{theorem} \label{thm:mse-lower-bound}
Suppose that $D \subseteq \R^d$ and $K \in \Sob(\nu)$ and $K_0 \in \Sob(\nu_0)$ for $\nu, \nu_0 > 0$.
For $x \in D$, let $\mathrm{dist}_\mathcal{X}(x) = \min_{x' \in \mathcal{X}} \lVert x - x' \rVert$ be the distance between $x$ and observation locations $\mathcal{X} \subseteq D$.
Then there is a constant $c > 0$, which does not depend on $x$ or $\mathcal{X}$, such that 
\begin{equation*}
        \MSE_\mathcal{X}(x) \geq c \cdot \mathrm{dist}_\mathcal{X}(x)^{2\nu} \quad \text{ and } \quad \MSE_\mathcal{X}^*(x) \geq c \cdot \mathrm{dist}_\mathcal{X}(x)^{2\nu_0}.
\end{equation*}
\end{theorem}
\begin{proof}
The proof is standard in the theory of approximation in Sobolev spaces; see \cite[Sec.\@~1.3.11]{Novak1988} or \citep{DeMarchiSchaback2010}.
We include it because the exact statement we need has proved difficult to locate.
Fix $x \in D$ and let $\delta = \mathrm{dist}_\mathcal{X}(x) = \min_{x^* \in \mathcal{X}} \lVert x - x^* \rVert$.
Define the bump function $\phi$ as
\begin{equation*}
        \phi(y) = \exp\bigg(\! -\frac{1}{1 - \lVert y \rVert^2} \bigg) \:\: \text{ if } \:\: \lVert y \rVert < 1 \quad \text{ and } \quad \phi(y) = 0 \:\: \text{ if } \:\: \lVert y \rVert \geq 1 
\end{equation*}
and $g$ as $g(y) = \phi((y-x)/\delta)$.
Since $\phi$ is supported on the unit ball, $\phi(0) = 1$, and $\lVert x - x' \rVert \geq \delta$ for every $x' \in \mathcal{X}$, the function $g$ vanishes at $x'$ for all $x' \in \mathcal{X}$ and takes value one at $y = x$.
Let
\begin{equation*}
        e_\mathcal{X}(x; H^\alpha(D)) = \sup_{ \lVert f \rVert_{H^\alpha(D)} \leq 1} \lvert f(x) - (I_\mathcal{X} f)(x) \rvert
\end{equation*}
be the worst-case error in a Sobolev space of order $\alpha > d/2$.
Here $I_\mathcal{X} f$ is the kernel interpolant in~\eqref{eq:kernel-interpolant} constructed using the observation locations $\mathcal{X}$.
Being infinitely differentiable and compactly supported, the bump function $\phi$ is an element of every Sobolev space.
Since $g$ is obtained via scaling and translation of $\phi$, it too is an element of every Sobolev space.
Therefore $g / \lVert g \rVert_{H^\alpha(D)}$ is an element of $H^\alpha(D)$ with unit norm.
Moreover, $I_\mathcal{X} g \equiv 0$ because $g$ vanishes on $\mathcal{X}$ and, by~\eqref{eq:kernel-interpolant}, the kernel interpolant is a linear combination of $g(x')$ for $x' \in \mathcal{X}$.
Thus
\begin{equation*}
        e_\mathcal{X}(x; H^\alpha(D))^2 \geq \frac{\lvert g(x) \rvert^2}{\lVert g \rVert_{H^\alpha(D)}^2} = \frac{1}{\lVert g \rVert_{H^\alpha(D)}^2} .
\end{equation*}
We are left to estimate the Sobolev norm of $g$. 
Observe that $\lvert (\mathcal{F} g)(\omega) \rvert^2 = \delta^{2d} \lvert (\mathcal{F} \phi)(\delta \omega) \rvert^2$.
When $\delta \leq 1$, Equation~\eqref{eq:Bessel-norm} and a change of variables give
\begin{equation*}
        \begin{split}
        \lVert g \rVert_{H^\alpha(D)}^2 \leq \lVert g \rVert_{H^\alpha(\R^d)}^2        
        &= \frac{\delta^{2d}}{(2\pi)^{d/2}}\int_{\R^d} \lvert (\mathcal{F} \phi)(\delta \omega) \rvert^2 (1 + \lVert \omega \rVert^2)^\alpha \dif \omega \\
        &= \frac{\delta^{d}}{(2\pi)^{d/2}}\int_{\R^d} \lvert (\mathcal{F} \phi)(\omega) \rvert^2 \bigg(1 + \frac{\lVert \omega \rVert^2}{\delta^2} \bigg)^\alpha \dif \omega \\
        &\leq \frac{\delta^{d-2\alpha}}{(2\pi)^{d/2}}\int_{\R^d} \lvert (\mathcal{F} \phi)(\omega) \rvert^2 (1 + \lVert \omega \rVert^2 )^\alpha \dif \omega \\
        &= \delta^{d - 2\alpha} \lVert \phi \rVert_{H^\alpha(\R^d)}^2 .
        \end{split}
\end{equation*}
Consequently, \smash{$e_\mathcal{X}(x; H_2^\alpha(D))^2 \geq \lVert \phi \rVert_{H^\alpha(\R^d)}^{-2} \delta^{2\alpha - d}$}.
Recall~\eqref{eq:wce-mse-1} and~\eqref{eq:wce-mse-2}, which are clearly valid for a general set of observation locations, $\mathcal{X}$.
Because $K$ and $K_0$ are Sobolev kernels of orders $\nu$ and $\nu_0$, the claim follows by plugging $\alpha = \nu + d/2$ and $\alpha = \nu_0 + d/2$ in the preceding equation and using norm-equivalence.
\end{proof}

The proofs of our main results rely on the following corollary.
As in \Cref{sec:scale-estimation}, we use the subscript $n \setminus k$ to indicate that the $k$th point has been removed from the dataset.

\begin{corollary} \label{cor:MSE-rates}
Suppose that $D$ is a bounded open convex set and $K \in \Sob(\nu)$ and $K_0 \in \Sob(\nu_0)$ for $\nu \geq \nu_0 > 0$.
If the sequence $(x_i)_{i=1}^\infty \subset D$ is quasi-uniform, then
\begin{equation*}
        \MSE_{n \setminus k}(x_k) \asymp \MSE_{n-1}(x_n) \asymp n^{-2\nu/d} \quad \text{and} \quad \MSE_{n \setminus k}^*(x_k) \asymp \MSE_{n-1}^*(x_n) \asymp n^{-2\nu_0/d}
\end{equation*}
for every $k$.
The hidden constants do not depend on $k$.
\end{corollary}
\begin{proof}
Quasi-uniformity of a sequence is not affected by the removal of one point when $D$ is convex.
The upper bounds thus follow immediately from \Cref{thm:mse-upper-bound}.
The lower bounds follow by selecting $\mathcal{X} = \{x_i\}_{i=1}^n \setminus \{x_k\}$ and $\mathcal{X} = \{x_i\}_{i=1}^{n-1}$ in \Cref{thm:mse-lower-bound} and noting that \smash{$\mathrm{dist}_\mathcal{X}(x_k) \geq 2 q_n \gtrsim n^{-1/d}$} by the definitions of separation radius and quasi-uniformity. 
As the constant $c$ in \Cref{thm:mse-lower-bound} does not depend on $x$ or $\mathcal{X}$, all constants are independent of $k$.
\end{proof}

The following theorems concern mean-square errors for periodic Sobolev kernels.
The periodic Sobolev kernel of order $\alpha > 1/2$ is defined in~\eqref{eq:korobov-kernel-alpha}. 
Its RKHS is the periodic Sobolev space of order $\alpha$, which has the characterisation
\begin{equation*} 
H(K) = H_{\textup{per}}^{\alpha}([0,1]) = \Set[\bigg]{ f = \sum_{k \in \Z} \hat{f}(k) \varphi_k }{ \lVert f \rVert_{H_{\textup{per}}^\alpha([0,1])}^2 = \lvert \hat{f}(0) \rvert^2 + \sum_{k \neq 0} \lvert k \rvert^{2\alpha} \lvert \hat{f}(k) \rvert^2 < \infty }
\end{equation*}
in terms of the rate of decrease of the Fourier coefficients $\hat{f}(k) = \int_0^1 f(x) e^{-2\pi \mathrm{i} k x} \dif x$.
For $\alpha \in \N$, the RKHS consists of those functions in the Sobolev space $H^\alpha([0,1])$ whose derivatives up to order $\alpha-1$ are periodic and its norm is equivalent to the classical Sobolev norm in~\eqref{eq:sobolev-norm}.

\begin{theorem} \label{thm:periodic-lower-bounds}
Suppose that $D = [0, 1]$ and that $K$ and $K_0$ are periodic Sobolev kernels of orders $\alpha, \alpha_0 > 1/2$.
For $x \in D$, let $\mathrm{dist}_\mathcal{X'}(x) = \min_{x' \in \mathcal{X}'} \lvert x - x' \rvert$ be the distance between $x$ and the set $\mathcal{X}' = \mathcal{X} \cup \{0, 1\} \subseteq D$.
Then there is a constant $c > 0$, which does not depend on $x$ or $\mathcal{X}$, such that 
\begin{equation*}
        \MSE_\mathcal{X}(x) \geq c \cdot \mathrm{dist}_\mathcal{X'}(x)^{2\alpha - 1} \quad \text{ and } \quad \MSE_\mathcal{X}^*(x) \geq c \cdot \mathrm{dist}_\mathcal{X'}(x)^{2\alpha_0 - 1}.
\end{equation*}
\end{theorem}
\begin{proof}
The proof is more or less identical to that of \Cref{thm:mse-lower-bound}.
Let $x \in (0, 1)$ and $\delta = \mathrm{dist}_{\mathcal{X}'}(x)$.
We can construct an infinitely differentiable function $\phi$ that vanishes outside of  $(-1, 1)$ and satisfies $\phi(0) = 1$.
Then the function $g(y) = \phi((y-x)/\delta)$ and all its derivatives vanish outside of $(x-\delta, x+\delta)$.
In particular, $g(x') = 0$ for every $x' \in \mathcal{X'}$.
Because $\mathcal{X}'$ includes the boundary, $g$ and all its derivatives are periodic.
Therefore \smash{$g \in H_{\textup{per}}^{\beta}([0, 1])$} for every $\beta \in \N$. 
Since \smash{$H_{\textup{per}}^{\beta}([0, 1]) \subset H_{\textup{per}}^{\gamma}([0, 1])$} if $\beta > \gamma$, the function $g$ is an element of \smash{$H(K) = H_{\textup{per}}^{\alpha}([0, 1])$}.
Because $g$ is in particular supported on $[0, 1]$, we have $\hat{g}(k) = (\mathcal{F} g)(k)$.
Thus
\begin{equation*}
        \begin{split}
        \lVert g \rVert_{H_{\textup{per}}^\alpha([0, 1])}^2 = \lvert \hat{g}(0) \rvert^2 + \sum_{k \neq 0} \lvert k \rvert^{2\alpha} \lvert \hat{g}(k) \rvert^2 
        &= \delta^2 \bigg( \int_\R \phi(x) \dif x \bigg)^2 + \delta^2 \sum_{k \neq 0} \lvert k \rvert^{2\alpha} \lvert (\mathcal{F} \varphi)(\delta k) \rvert^2.
        \end{split}
\end{equation*}
As $g \in H^\beta(\R)$ for every $\beta > 1/2$, there is $C > 0$ such that $\lvert (\mathcal{F} g)(\omega) \rvert^2 \leq C(1 + \lvert \omega \rvert^2 )^{-\alpha - 1}$ for all $\omega \in \R$.
With this bound it is straightforward to obtain the estimate \smash{$\lVert g \rVert_{H_{\textup{per}}^\alpha([0, 1])}^2 = O(\delta^{1-2\alpha})$} as $\delta \to 0$.
The rest of the proof follows that of \Cref{thm:mse-lower-bound}.
\end{proof}

\begin{theorem} \label{thm:periodic-upper-lower}
Suppose that $D = [0, 1]$ and that $K$ and $K_0$ are periodic Sobolev kernels of orders $\alpha \in \N$ and $\alpha_0 \in \R$ such that $2\alpha \geq \alpha_0 \geq \alpha > 1/2$.
If the sequence $(x_i)_{i=1}^\infty \subset D$ is quasi-uniform, then
\begin{equation*}
        \MSE_{n \setminus k}(x_k) \asymp \MSE_{n-1}(x_n) \asymp n^{-2\alpha + 1} \quad \text{and} \quad \MSE_{n \setminus k}^*(x_k) \asymp \MSE_{n-1}^*(x_n) \asymp n^{-2\alpha_0 + 1}
\end{equation*}
for every $k$. The hidden constants do not depend on $k$.
\end{theorem}
\begin{proof}
The lower bounds follow from \Cref{thm:periodic-lower-bounds} and quasi-uniformity.
Observe from~\eqref{eq:wce-explicit} that the squared worst-case error has the form
\begin{equation*}
        e_n(x; R)^2 = R(x, x) - 2 \b{r}_n(x) \b{c}_n(x) + \b{c}_n(x)^\T \b{R}_n \b{c}_n(x)
\end{equation*}
for a certain vector $\b{c}_n(x)$.
The right-hand side is a quadratic form that is minimised by setting $\b{c}_n(x) = \b{R}_n(x)^{-1} \b{r}_n(x)$, where $\b{r}_n(x) = (R(x, x_1), \ldots, R(x, x_n)) \in \R^n$.
This corresponds to setting $R = K$, so that the smallest possible worst-case error is obtained by constructing the interpolant using the reproducing kernel.
Because \smash{$H(K) = H_{\textup{per}}^\alpha([0, 1]) \subset H^\alpha([0,1])$} with norm-equivalence, the claimed upper bounds on $\MSE_{n \setminus k}(x_k) = e_{n \setminus k}(x_k; K)^2$ and $\MSE_{n-1}(x_n) = e_{n-1}(x_n; K)^2$ follow from~\eqref{eq:wce-mse-2} and \Cref{thm:arcangeli} ($\nu = \alpha - 1/2$ and $p = \infty$).

The upper bounds on the true mean-square errors use results from~\cite{KarvonenSantinWenzel2025}.
The definition of the worst-case error in~\eqref{eq:wce-def} implies that $\lvert g(x) - (I_n g)(x) \rvert \leq \lVert g \rVert_{H(K)} \, e_n(x; K)$
for all $g \in H(K)$.
Since $I_n(f - I_n f) \equiv 0$, selecting $g = f - I_n f$ yields
\begin{equation} \label{eq:superconvergence-f-Inf}
        \lvert f(x) - (I_n f)(x) \rvert \leq \lVert f - I_n f \rVert_{H(K)} \, e_n(x; K) 
\end{equation}
for all $f \in H(K)$.
Let $\theta \geq 0$.
The so-called \emph{$\theta$th power}, $H_\theta(K)$, of an RKHS $H(K)$ is obtained by raising its Mercer eigenvalues to power $\theta$~\cite[Sec.\@~4]{SteinwartScovel2012}.
By Example~18 in~\cite{KarvonenSantinWenzel2025}, for $H(K) = H_\textup{per}^\alpha([0, 1])$ the power space is conveniently simply \smash{$H_\theta(K) = H_\textup{per}^{\theta\alpha}([0, 1])$}.
If $f$ is an element of $H_\theta(K)$ for $\theta \in [1, 2]$, it follows from $H(K) \subset H^\alpha([0, 1])$, \Cref{thm:arcangeli} ($p = 2$ and $\nu = \alpha - 1/2$), and
Corollary~15 in~\cite{KarvonenSantinWenzel2025} that $\lVert f - I_n f \rVert_{H(K)} \leq C n^{-(\theta-1) \alpha} \lVert f \rVert_{H_\theta(K)}$ for all $n$ and a positive constant $C$ that does not depend on $f$.
By combining this with~\eqref{eq:superconvergence-f-Inf} and \Cref{thm:arcangeli} ($p = \infty$ and $\nu = \alpha - 1/2$) we obtain
\begin{equation*}
        \sup_{\lVert f \rVert_{H_\theta(K)} \leq 1} \lvert f(x) - (I_n f)(x) \rvert \lesssim n^{-\alpha + 1/2} \cdot n^{-(\theta-1) \alpha} = n^{-\theta \alpha + 1/2} 
\end{equation*}
for $\theta \in [1, 2]$.
The claimed upper bounds on the true mean-square errors are now obtained by setting $\alpha_0 = \theta \alpha$, so that \smash{$H_\theta(K) = H_\textup{per}^{\alpha_0}([0, 1])$}, and using~\eqref{eq:wce-mse-1}.
\end{proof}

\begin{theorem} \label{thm:MSE-bounds-periodic}
        Let $p \in (0, \infty)$.
        Suppose that $D = [0, 1]$ and that $K$ and $K_0$ are periodic Sobolev kernel of orders $\alpha \in \N$ and $\alpha_0 \in \R$ such that $2\alpha \geq \alpha_0 \geq \alpha > 1/2$.
        If the sequence $(x_i)_{i=1}^\infty \subset D$ is quasi-uniform, then
        \begin{equation*}
                \lVert \MSE_{n} \rVert_{L^p(D)} \asymp n^{-2\alpha + 1} \quad \text{ and } \quad \lVert \MSE^*_{n} \rVert_{L^p(D)} \asymp n^{-2\alpha_0 + 1} .
        \end{equation*}
\end{theorem}
\begin{proof}
The claim follows from \Cref{thm:periodic-lower-bounds,thm:periodic-upper-lower} and an argument similar to that used in the proof of \Cref{thm:MSE-rates-informal} in \Cref{sec:proofs-setting}.
\end{proof}

\subsection{Proofs for \Cref{sec:setting}} \label{sec:proofs-setting}

This section contains the proofs for \Cref{sec:setting}.

\begin{proof}[Proof of \Cref{thm:MSE-rates-informal}]
\Cref{thm:mse-upper-bound} and $D$ being bounded yield the upper bounds.
For the lower bounds we use \Cref{thm:mse-lower-bound}.
Consider $\MSE_n$ (the proof for $\MSE_n^*$ is identical).
Let $A_n = \Set{x \in D}{ \mathrm{dist}_{\R^d \setminus D}(x) \leq q_n}$
be the ``annulus'' of width $q_n$ \emph{inside} $D$.
Open balls $B(x_i, q_n)$ of radius $q_n$ centered at $x_i$ for $i = 1, \ldots n$ are disjoint by the definition of the separation radius, $q_n$. 
Therefore
\begin{equation*}        
        \int_D \MSE_n(x)^p \dif x \geq \sum_{x_i \in D \setminus A_n} \int_{B(x_i, q_n)} \MSE_n(x)^p \dif x ,
\end{equation*}
where summation is over the first $n$ observation locations and $x_i \in D \setminus A_n$ ensures that each ball is contained in $D$.
By \Cref{thm:mse-lower-bound} and the definition of the separation radius,
\begin{equation*}
        \int_{B(x_i, q_n)} \MSE_n(x)^p \dif x \geq c^p \int_{B(x_i, q_n)} \lVert x - x_i \rVert^{2\nu p} \dif x = c^p \int_{B(0, q_n)} \lVert x \rVert^{2\nu p } \dif x .
\end{equation*}
The integral is
\begin{equation*}
        \int_{B(0, q_n)} \lVert x \rVert^{2\nu p } \dif x = C_d \int_0^{q_n} r^{2\nu p} \cdot r^{d-1} \dif r = \frac{C_d}{2\nu p  + d} \, q_n^{2\nu p + d} 
\end{equation*}
for a constant $C_d > 0$ depending only on $d$.
Since $q_n \asymp n^{-1/d}$ by quasi-uniformity, we obtain
\begin{equation*}
        \int_D \MSE_n(x)^p \dif x \geq \frac{C_d \, c^p }{2\nu p  + d} \sum_{x_i \in D \setminus A_n} q_n^{2\nu p + d} \gtrsim n^{-2\nu p / d - 1}\sum_{x_i \in D \setminus A_n} 1 .
\end{equation*}
We now need to estimate how many of the first $n$ observation locations are in $D \setminus A_n$.
Because $D$ is convex, there are $x_0 \in D$ and $\delta > 0$ such that $B(x_0, \delta) \subset D \setminus A_n$ for all sufficiently large~$n$. 
By quasi-uniformity, the number of observation locations in $B(x_0, \delta)$ must grow as $n$ (see~\cite[Sec.\@~14.1]{Wendland2005} for such arguments).
Thus
\begin{equation*}
        \int_D \MSE_n(x)^p \dif x \geq n^{-2\nu p / d - 1}\sum_{x_i \in D \setminus A_n} 1 \geq n^{-2\nu p / d - 1}\sum_{x_i \in B(x_0, \delta )} 1 \gtrsim n^{-2\nu p / d} ,
\end{equation*} 
from which it follows that
\begin{equation*}
        \lVert \MSE_n \rVert_{L^p(D)} = \bigg( \int_D \MSE_n(x)^p \dif x \bigg)^{1/p} \gtrsim n^{-2\nu / d}. \qedhere
\end{equation*}
\end{proof}

\begin{proof}[Proof of \Cref{prop:rate-unbiased-lambda}]
The claim follows from~\eqref{eq:wce-mse-1} and~\eqref{eq:wce-mse-2} and the fact that the RKHS of a Matérn kernel is norm-equivalent to a Sobolev space with norm-equivalence constants that are bounded away from zero and infinity when the correlation length varies on a bounded interval~\cite[Lem.\@~3.4]{Teckentrup2020}.
This is easy to verify by comparing~\eqref{eq:matern-FT} and~\eqref{eq:Bessel-norm}.
\end{proof}

\subsection{Proofs for \Cref{sec:results}} \label{sec:proofs-results-sigma}

This section contains the proofs for \Cref{sec:results}.
The following lemma allows transforming rates in expectation to rates in probability.

\begin{lemma} \label{lemma:Prob-E-lemma}
Let $(a_n)_{n=1}^\infty$ be a non-negative sequence and let $\hat{\sigma}_n^2$ stand for either $\hat{\sigma}_{\ML,n}^2$ or \smash{$\hat{\sigma}_{\CV, n}^2$}. Then \smash{$\hat{\sigma}_n^2 \asymp_{\mathbb{P}} a_n$} if \smash{$\mathbb{E}[ \hat{\sigma}_n^2 ] \asymp a_n$}.
\end{lemma}
\begin{proof}
        Markov's inequality and the assumption \smash{$\mathbb{E}[ \hat{\sigma}_n^2 ] \asymp a_n$} yield
        \begin{align*}
        \limsup_{n \to \infty}
        \mathbb{P} \big(
        \hat{\sigma}_{n}^2 \ge
        \varepsilon^{-1} a_n
        \big)
        \le 
        \varepsilon \cdot
        \limsup_{n \to \infty}
        \mathbb{E}[ \hat{\sigma}_{n}^2 ] \cdot a_n^{-1} \to 0
        \quad \text{ as } \quad \varepsilon \to 0 .
        \end{align*}
        We are left to prove that $\limsup_{n \to \infty} \mathbb{P}( \hat{\sigma}_{n}^2 \leq \varepsilon a_n) \to 0$ as $\varepsilon \to 0$.

        Consider first the maximum likelihood estimator and denote \smash{$\hat{\sigma}_n = \hat{\sigma}_{\ML,n}$} for brevity. 
        Let $\b{K}_{0,n}^{\smash{1/2}}$ be the unique positive-definite matrix square root of $\b{K}_{0,n}$, the covariance of the observations $\b{y}_n$. 
        Perform the diagonalisation
        \begin{equation*}
                \frac{1}{n}\b{K}_{0,n}^{1/2} \b{K}_n^{-1} \b{K}_{0,n}^{1/2} = \b{P}_n^\T \b{D}_n \b{P}_n,
        \end{equation*}
        where $\b{P}_n$ is orthonormal and $\b{D}_n = \mathrm{diag}(\lambda_{n,1}, \ldots, \lambda_{n,n})$ is diagonal with non-increasing and non-negative diagonal elements.
        Let $\b{K}_{0,n}^{-\smash{1/2}}$ be the inverse of $\b{K}_{0,n}^{\smash{1/2}}$.
        Then~\eqref{eq:sigma-ML} yields
        \begin{equation} \label{eq:diagonalization}
            \hat{\sigma}_{n}^2 = 
            \big(
        \b{P}_n
             \b{K}_{0,n}^{-1/2} 
        \b{y}_n
            \big)^\T
        \b{D}_n
          \big(
        \b{P}_n 
          \b{K}_{0,n}^{-1/2}
        \b{y}_n
            \big) 
            =
        \sum_{i=1}^n
        \lambda_{n,i}
        z_{n,i}^2,
        \end{equation}
        where $z_{n,1}, \ldots , z_{n,n}$ are independent standard normal. 
        From \smash{$\mathbb{E}[ \hat{\sigma}_{n}^2 ] \asymp a_n$} it follows that
        \begin{equation} \label{eq:sum:lambda}
        \sum_{i=1}^n 
        \lambda_{n,i} 
        \asymp  a_n.
        \end{equation}
        Let $\varepsilon > 0$.
        Suppose first that $\lambda_{n,1} \ge \varepsilon \, a_n$.
        Then
        \begin{equation} \label{eq:lambda-estimate-1}
        \mathbb{P}
        \big(  
          \hat{\sigma}_{n}^2 
          \le 
          \varepsilon^2 
          a_n 
        \big)
        \le 
        \mathbb{P}
        \big(  
         \varepsilon \, a_n
         z_{n,1}^2
          \le 
          \varepsilon^2 
          a_n
        \big)
        = 
        \mathbb{P}\big( z_{n,1}^2 \leq \varepsilon \big) 
        =
        \int_{- \sqrt{\varepsilon}}^{\sqrt{\varepsilon}} \,
        g(t)
        \dif t,
        \end{equation}
        with $g$ the standard Gaussian density function. 
        Suppose then that $\lambda_{n,1} < \varepsilon \, a_n$.
        From \smash{$\mathbb{E}[ \hat{\sigma}_{n}^2 ] \asymp a_n$} it follows that there is a constant $c > 0$ such that
        \begin{equation*}
                \varepsilon^2 a_n - \mathbb{E}[ \hat{\sigma}_{n}^2 ] \leq -\frac{1}{2} \mathbb{E}[ \hat{\sigma}_{n}^2 ]\leq - c \, a_n 
        \end{equation*}
        for any sufficiently small $\varepsilon$ and all $n$.
        Therefore Chebyshev's inequality and~\eqref{eq:diagonalization} yield     
        \begin{align*}
                \mathbb{P}\big( \hat{\sigma}_{n}^2 \leq \varepsilon^2 a_n \big) 
                &= \mathbb{P}\big( \hat{\sigma}_{n}^2 - \mathbb{E}[ \hat{\sigma}_{n}^2 ] \leq \varepsilon^2 a_n - \mathbb{E}[ \hat{\sigma}_{n}^2 ] \big) \\
                &\hspace{2cm}\leq  \mathbb{P}\big( \hat{\sigma}_{n}^2 - \mathbb{E}[ \hat{\sigma}_{n}^2 ] \leq - c \, a_n \big) 
                \leq \frac{\Var[\hat{\sigma}_{n}^2]}{c^2 a_n^2} 
                = \frac{2\sum_{i=1}^n \lambda_{n,i}^2}{c^2 a_n^2} .
        \end{align*}
        From $\lambda_{n,1} \geq \lambda_{n,2} \geq \cdots \geq \lambda_{n,n}$ and $\lambda_{n,1} < \varepsilon \, a_n$ we obtain
        \begin{equation} \label{eq:lambda-estimate-2}
                \mathbb{P}\big( \hat{\sigma}_{n}^2 \leq \varepsilon^2 a_n \big) \leq \frac{2\sum_{i=1}^n \lambda_{n,i}^2}{c^2 a_n^2} \leq \frac{2 \lambda_{n,1}\sum_{i=1}^n \lambda_{n,i}}{c^2 a_n^2} < \frac{2 \varepsilon \sum_{i=1}^n \lambda_{n,i}}{c^2 a_n} .
        \end{equation}
        Combining the estimates~\eqref{eq:lambda-estimate-1} and~\eqref{eq:lambda-estimate-2} for the cases $\lambda_{n,1} \geq \varepsilon a_n$ and $\lambda_{n,1} < \varepsilon a_n$ gives
        \begin{equation*}
               \mathbb{P}\big( \hat{\sigma}_{n}^2 \leq \varepsilon^2 a_n \big) \leq \int_{- \sqrt{\varepsilon}}^{\sqrt{\varepsilon}} \, g(t) \dif t +  \frac{2 \varepsilon \sum_{i=1}^n \lambda_{n,i}}{c^2 a_n}
        \end{equation*}
        for any sufficiently small $\varepsilon$ and all $n$.
        From~\eqref{eq:sum:lambda} it follows that $\limsup_{n \to \infty} \mathbb{P}( \hat{\sigma}_{n}^2 \leq \varepsilon a_n) \to 0$ as $\varepsilon \to 0$, which concludes the proof for the maximum likelihood estimator.

        The proof for cross-validation is identical except for the use of~\eqref{eq:cv-alternate} and the diagonalisation
        \begin{equation*}
                \frac{1}{n}\b{K}_{0,n}^{1/2} \b{K}_n^{-1} (\diag \b{K}_n^{-1})^{-1 }\b{K}_n^{-1} \b{K}_{0,n}^{1/2} = \b{P}_n^\T \b{D}_n \b{P}_n . \qedhere
        \end{equation*}
\end{proof}

\begin{proof}[Proof of \Cref{thm:estimator-asymptotics}]
We prove~\eqref{eq:sigma-E-asymp} for expectations.
The probabilistic asymptotics then follow from \Cref{lemma:Prob-E-lemma}.
Consider first the maximum likelihood estimator
\begin{equation*} 
        \hat{\sigma}_{\ML,n}^2 = \frac{1}{n} \sum_{k=1}^n \frac{( y_k - \mu_{k-1}(x_k))^2}{V_{k-1}(x_k)} .
\end{equation*}
Because $\b{y}_n \sim \mathrm{N}(0, \b{K}_{0,n})$, taking expectation gives
\begin{equation} \label{eq:sigma-ML-E-app}
        \E[ \hat{\sigma}_{\ML,n}^2 ] = \frac{1}{n} \sum_{k=1}^n \frac{\E[ y_k - \mu_{k-1}(x_k)]^2}{V_{k-1}(x_k)} = \frac{1}{n} \sum_{k=1}^n \frac{\MSE_{k-1}^*(x_k)}{\MSE_{k-1}(x_k)},
\end{equation}
where we used~\eqref{eq:MSE} and~\eqref{eq:MSE-true}.
\Cref{cor:MSE-rates} yields
\begin{equation*}
        \E[ \hat{\sigma}_{\ML,n}^2 ] = \frac{1}{n} \sum_{k=1}^n \frac{\MSE_{k-1}^*(x_k)}{\MSE_{k-1}(x_k)} \asymp \frac{1}{n} \sum_{k=1}^n \frac{k^{-2\nu_0/d}}{k^{-2\nu/d}} = \frac{1}{n} \sum_{k=1}^n k^{2(\nu-\nu_0)/d}.
\end{equation*}
The asymptotics $\sum_{k=1}^n k^{2(\nu-\nu_0)/d} \asymp n^{2(\nu-\nu_0)/d + 1}$ implies~\eqref{eq:sigma-E-asymp} for the maximum likelihood estimator.
Consider then cross-validation.
The expectation of the scale estimator is
\begin{equation} \label{eq:sigma-CV-E-app}
        \E[\sigma_{\CV,n}^2] = \E\bigg[ \frac{1}{n} \sum_{k=1}^n \frac{( y_k - \mu_{n \setminus k}(x_k))^2}{V_{n \setminus k}(x_k)} \bigg] = \frac{1}{n} \sum_{k=1}^n \frac{\MSE_{n \setminus k}^*(x_k)}{\MSE_{n \setminus k}(x_k)} .
\end{equation}
\Cref{cor:MSE-rates} then yields~\eqref{eq:sigma-E-asymp} for cross-validation:
\begin{equation*}
        \E[ \hat{\sigma}_{\CV,n}^2 ] = \frac{1}{n} \sum_{k=1}^n \frac{\MSE_{n \setminus k}^*(x_k)}{\MSE_{n \setminus k}(x_k)} \asymp \frac{1}{n} \sum_{k=1}^n \frac{n^{-2\nu_0/d}}{n^{-2\nu/d}} = \frac{1}{n} \sum_{k=1}^n n^{2(\nu-\nu_0)/d} = n^{2(\nu-\nu_0)/d} . \qedhere
\end{equation*}
\end{proof}

\begin{proof}[Proof of \Cref{thm:mle-undersmooth-exp}]
We prove parts (a) and (b) of the theorem separately.

\emph{Part (a).}
The identity $\E[ \sigma_{\ML,n}^2 ] = n^{-1} \tr( \b{K}_n^{-1} \b{K}_{0,n} )$ was obtained in~\eqref{eq:sigma-ML-unbiased}. 
We are left to prove that $H(K_0) \subseteq H^\alpha(D)$ for $\alpha > \nu + d$ implies that the trace tends to a finite positive constant.
From~\eqref{eq:sigma-ML-E-app}, together with~\eqref{eq:wce-def} and~\eqref{eq:wce-mse-1}, it follows that
\begin{equation} \label{eq:trace-sup-explicit}
        \tr( \b{K}_n^{-1} \b{K}_{0,n} ) = \sum_{k=1}^n \frac{\MSE_{k-1}^*(x_k)}{\MSE_{k-1}(x_k)} = \sum_{k=1}^n \frac{\sup_{\lVert f \rVert_{H(K_0)} \leq 1} \lvert f(x_k) - (I_{k-1} f)(x_k) \rvert^2}{\sup_{\lVert f \rVert_{H(K)} \leq 1} \lvert f(x_k) - (I_{k-1} f)(x_k) \rvert^2} .
\end{equation}
Therefore the trace is positive and non-decreasing.
Since $H^\alpha(D)$ is an RKHS (recall \Cref{rmk:equivalence-classes}), $H(K_0)$ is continuously embedded in $H^\alpha(D)$ by Theorem~5.1 in~\citep{Paulsen2016}.
As each numerator in~\eqref{eq:trace-sup-explicit} is a supremum over the unit ball of $H(K_0)$ and we wish to bound the trace from above, we may thus assume that $H(K_0)$ is norm-equivalent to $H^\alpha(D)$.
In fact, again by norm-equivalence, we may assume that $K_0$ is continuous on $D \times D$.
Recall that $K$ was assumed continuous on $D \times D$.
Suppose for a moment that $(x_i)_{i=1}^\infty$ is dense in $D$.
Under these continuity and denseness assumptions we can invoke Proposition~4.4 in \citep{LukicBeder2001}, which states that $\lim_{n \to \infty} \tr( \b{K}_n^{-1} \b{K}_{0,n} ) = \tr(L)$, where $L$ is a certain dominance operator between the RKHSs of $K$ and $K_0$.
By Theorem~7.4 in \citep{LukicBeder2001}, the finiteness of $\tr(L)$ is equivalent to the existence of a version of $X_0 \sim \GP(0, K_0)$ whose sample paths are almost surely in $H(K)$.
In the case of Sobolev spaces, Corollary~5.7 in \citep{Steinwart2019} states that a Gaussian process $X_0$ with covariance $K_0 \in \Sob(\nu_0)$ has a version whose samples are almost surely in \smash{$H_2^\beta(D)$} if and only if $\nu_0 > \beta$, which is equivalent to $\alpha > \nu + d$.
Thus the claim holds when the sequence of observation locations is dense in $D$.
If the sequence is not dense, for each $n$ we can concatenate $(x_i)_{i=1}^n$ to the beginning of some dense sequence, so that the resulting trace will tend to $\tr(L)$.
The trace being non-decreasing, we are assured that $\tr( \b{K}_n^{-1} \b{K}_{0,n} ) \leq \tr(L) < \infty$ for every $n$.

\emph{Part (b).} Under these assumptions almost all sample paths of $X_0$ are elements of the Sobolev space $H^{\nu_0-\varepsilon}(D)$ for every $\varepsilon > 0$.
For this relatively well known result, see Corollary~5.7 in \citep{Steinwart2019} and Section~4.4 in \citep{Kanagawa2018} (see also \cite{Driscoll1973, Henderson2024, Karvonen2023b, LukicBeder2001, Scheuerer2010}).
By the continuity assumption we do not have to consider a version of the stochastic process.
Since $K \in \Sob(\nu)$, its RKHS is norm-equivalent to \smash{$H^{\nu+d/2}(D)$} and thus contains almost all samples by the assumption $\nu_0 > \nu + d/2$.
We may thus apply Proposition~3.1 in \citep{Karvonen2020} to almost every sample path.
This proposition applies to functions which do not vanish everywhere. 
Sample paths satisfy this because $X_0(x) \sim \mathrm{N}(0, K_0(x, x))$ for every $x \in D$ and $K_0(x, x) > 0$ by positive-definiteness.
\end{proof}

It seems probable that part~(a) of \Cref{thm:mle-undersmooth-exp} has a less tortuous proof that does not require invoking results on sample path properties of Gaussian processes. 

\begin{proof}[Proof of \Cref{thm:periodic}]
We proceed as in the proof of \Cref{thm:estimator-asymptotics}, except that we use \Cref{thm:periodic-upper-lower} to control the mean-square errors.
For the maximum likelihood estimator we get
\begin{equation*}
        \E[ \hat{\sigma}_{\ML,n}^2 ] = \frac{1}{n} \sum_{k=1}^n \frac{\MSE_{k-1}^*(x_k)}{\MSE_{k-1}(x_k)} \asymp \frac{1}{n} \sum_{k=1}^n \frac{k^{-2\alpha_0+1}}{k^{-2\alpha+1}} = \frac{1}{n} \sum_{k=1}^n k^{2(\alpha - \alpha_0)}.
\end{equation*}
The claim follows from the asymptotics
\begin{equation*}
        \sum_{k=1}^n k^{2(\alpha - \alpha_0)} \asymp 
        \begin{cases}
                n^{2(\alpha - \alpha_0) + 1} &\text{ if } \quad \alpha_0 < \alpha + 1/2 , \\
                \log n &\text{ if } \quad \alpha_0 = \alpha + 1/2, \\
                1 &\text{ if } \quad \alpha_0 > \alpha + 1/2 .
        \end{cases}
\end{equation*}
The second case uses the asymptotics $\sum_{k=1}^n k^{-1} \sim \log n$ for the $n$th harmonic number.
For cross-validation we get
\begin{equation*}
        \E[ \hat{\sigma}_{\CV,n}^2 ] = \frac{1}{n} \sum_{k=1}^n \frac{\MSE_{n \setminus k}^*(x_k)}{\MSE_{n \setminus k}(x_k)} \asymp \frac{1}{n} \sum_{k=1}^n \frac{n^{-2\alpha_0+1}}{n^{-2\alpha+1}} = \frac{1}{n} \sum_{k=1}^n n^{2(\alpha-\alpha_0)} = n^{2(\alpha - \alpha_0)} . \qedhere
\end{equation*}
\end{proof}

\subsection{Proofs for \Cref{sec:numerical}} \label{sec:proofs-numerical}

This section contains the proofs for \Cref{sec:numerical}.

\begin{proof}[Proof of \Cref{prop:smoothness:LS}]

Recall that $n_1 = n_1(n),\ldots,n_m = n_m(n)$.
It suffices to show that $\hat{\beta}_n \to \beta_0$ in probability,
where $\beta_0 = 2(\nu-\nu_0)/d$.
Let $\varepsilon>0$ be fixed and consider the event $\{ \hat{\beta}_n \ge \beta_0 + \varepsilon \}$.
Define 
\[
S_{n}( a,\beta )
=
\sum_{i=1}^m ( a + \beta \log n_i - \log \hat{\sigma}_{n_i}^2 )^2. 
\]
Since $m$ is fixed, \Cref{thm:estimator-asymptotics} implies that
\[
S_n(\hat{a}_n,\hat{\beta}_n)
\le 
S_n( 0,\beta_0 )
= 
\sum_{i=1}^m
( 
\beta_0 \log n_i 
- 
\log \hat{\sigma}_{n_i}^2
)^2
=
O_{\mathbb{P}}(1) .
\]
\Cref{thm:estimator-asymptotics} also implies that
\begin{align} \label{eq:sum:square}
S_n(\hat{a}_n,\hat{\beta}_n) 
&\ge
\big( 
\hat{a}_n
+ 
\hat{\beta}_n \log n_1 
- 
\log \hat{\sigma}_{n_1}^2
\big)^2
+
\big( 
\hat{a}_n
+ 
\hat{\beta}_n \log n_m 
- 
\log \hat{\sigma}_{n_m}^2
\big)^2
\notag 
\\ 
& = 
\big( 
\hat{a}_n
+ 
[ \hat{\beta}_n
-
\beta_0]
\log n_1 
+ O_{\mathbb{P}}(1)
\big)^2
+
\big( 
\hat{a}_n
+ 
[ \hat{\beta}_n
-
\beta_0]
\log n_m 
+ O_{\mathbb{P}}(1)
\big)^2.
\end{align}
Write $(t)_- = \max\{0,-t\}$ for $t \in \mathbb{R}$. 
If $-\hat{a}_n \geq (\hat{\beta}_n - \beta_0) (\log n_1 + \log n_m) / 2$
we obtain, still under the event \smash{$\{ \hat{\beta}_n \ge \beta_0 + \varepsilon \}$}, that
\begin{align*}
S_n(\hat{a}_n,\hat{\beta}_n) 
\ge &
\bigg( 
(\hat{\beta}_n
-
\beta_0)
\log n_1 
-
( \hat{\beta}_n
-
\beta_0)
\frac{  \log n_1 
+ 
 \log n_m 
}{2}
+ O_{\mathbb{P}}(1)
\bigg)_-^2
\\
= &
\bigg( 
- (\hat{\beta}_n
-
\beta_0)
\frac{\log n_m - \log n_1}{2} 
+ O_{\mathbb{P}}(1)
\bigg)_-^2
\\ 
\ge &
\bigg( 
- 
\varepsilon
\frac{\log n_m - \log n_1}{2} 
+ O_{\mathbb{P}}(1)
\bigg)_-^2. 
\end{align*}
Hence the event
\begin{equation} \label{eq:the:event}
   \bigg\{ 
\hat{\beta}_n \ge \beta_0 + \varepsilon \: \text{ and } \:
-\hat{a}_n
\ge 
(\hat{\beta}_n
-
\beta_0)
\frac{  \log n_1 
+ 
 \log n_m 
}{2}
   \bigg\}
\end{equation}
implies the event
\[
\bigg\{
\bigg( 
- 
\varepsilon
\frac{\log n_m - \log n_1}{2} 
+ O_{\mathbb{P}}(1)
\bigg)_-^2
\le 
O_{\mathbb{P}}(1)
\bigg\}.
\]
The latter event has vanishing probability as $n \to \infty$ since $\log n_m - \log n_1 \to \infty$ by assumption. Hence, the probability of the event in \eqref{eq:the:event} also vanishes. 

Write $(t)_+ = \max\{0,t\}$ for $t \in \mathbb{R}$. If
$
- \hat{a}_n
\le 
( \hat{\beta}_n
-
\beta_0)
(  \log n_1 
+ 
 \log n_m 
)/2, 
$
we now obtain from~\eqref{eq:sum:square}, also under the event $\{ \hat{\beta}_n \ge \beta_0 + \varepsilon \}$, that

\begin{align*}
S_n(\hat{a}_n,\hat{\beta}_n) 
\ge &
\bigg( 
( \hat{\beta}_n
-
\beta_0)
\log n_m 
-
( \hat{\beta}_n
-
\beta_0)
\frac{  \log n_1 
+ 
 \log n_m 
}{2}
+ O_{\mathbb{P}}(1)
\bigg)_+^2
\\
= &
\bigg( 
 ( \hat{\beta}_n
-
\beta_0)
\frac{\log n_m - \log n_1}{2} 
+ O_{\mathbb{P}}(1)
\bigg)_+^2 
\\ 
\ge &
\bigg( 
\varepsilon
\frac{\log n_m - \log n_1}{2} 
+ O_{\mathbb{P}}(1)
\bigg)_+^2. 
\end{align*}
Hence the event
\begin{equation} \label{eq:the:event:two}
   \bigg\{ 
\hat{\beta}_n \ge \beta_0 + \varepsilon,
- \hat{a}_n
\le 
( \hat{\beta}_n
-
\beta_0)
\frac{  \log n_1 
+ 
 \log n_m 
}{2}
   \bigg\}
\end{equation}
implies the event
\[
\bigg\{
\bigg( 
\varepsilon
\frac{\log n_m - \log n_1}{2} 
+ O_{\mathbb{P}}(1)
\bigg)_+^2
\le 
O_{\mathbb{P}}(1)
\bigg\}.
\]
The latter event has vanishing probability as $n \to \infty$ since $\log n_m - \log n_1 \to \infty$. Hence, the probability of the event in \eqref{eq:the:event:two} also vanishes. 
In the end we have shown
$
\mathbb{P}( \hat{\beta}_n \ge \beta_0 + \varepsilon)
\to
0
$.
We can show similarly
\smash{$
\mathbb{P}
( 
\hat{\beta}_n \le \beta_0 - \varepsilon
)
\to
0
$},
which concludes the proof.
\end{proof}

\begin{acks}[Acknowledgments]
We are grateful to Motonobu Kanagawa for numerous discussions and suggestions that helped to determine the direction of this research and to Anatoly Zhigljavsky and Luc Pronzato for a counter-example to a point-wise version of \Cref{thm:MSE-rates-informal}.
\end{acks}

\begin{funding}
TK was generously supported by the Research Council of Finland projects 338567 (``Scalable, adaptive and reliable probabilistic integration''), 359183 (``Flagship of Advanced Mathematics for Sensing, Imaging and Modelling''), and 368086 (``Inference and approximation under misspecification'').
FB was supported by the Project GAP (ANR-21-CE40-0007) of the French National Research Agency (ANR).
Part of the research was carried out during a research visit partially funded by the Maupertuis Program of the French Institute in Finland, the Embassy of France to Finland, and the Finnish Society of Sciences and Letters.

\end{funding}


\begin{thebibliography}{80}

\bibitem{AdamsFournier2003}
\begin{bbook}[author]
\bauthor{\bsnm{Adams},~\bfnm{R.~A.}\binits{R.~A.}} \AND \bauthor{\bsnm{Fournier},~\bfnm{J.~J.~F.}\binits{J.~J.~F.}}
(\byear{2003}).
\btitle{{S}obolev Spaces},
\bedition{2nd} ed.
\bpublisher{Academic Press}.
\end{bbook}
\endbibitem

\bibitem{AdlerTaylor2007}
\begin{bbook}[author]
\bauthor{\bsnm{Adler},~\bfnm{R.~J.}\binits{R.~J.}} \AND \bauthor{\bsnm{Taylor},~\bfnm{J.~E.}\binits{J.~E.}}
(\byear{2007}).
\btitle{Random Fields and Geometry}.
\bpublisher{Springer}.
\end{bbook}
\endbibitem

\bibitem{anderes2010consistent}
\begin{barticle}[author]
\bauthor{\bsnm{Anderes},~\bfnm{Ethan}\binits{E.}}
(\byear{2010}).
\btitle{On the consistent separation of scale and variance for {Gaussian} random fields}.
\bjournal{The Annals of Statistics}
\bvolume{38}
\bpages{870--893}.
\end{barticle}
\endbibitem

\bibitem{Arcangeli2007}
\begin{barticle}[author]
\bauthor{\bsnm{Arcang{\'e}li},~\bfnm{R.}\binits{R.}}, \bauthor{\bparticle{de} \bsnm{Silanes},~\bfnm{M.~C.~L.}\binits{M.~C.~L.}} \AND \bauthor{\bsnm{Torrens},~\bfnm{J.~J.}\binits{J.~J.}}
(\byear{2007}).
\btitle{An extension of a bound for functions in {S}obolev spaces, with applications to $(m,s)$-spline interpolation and smoothing}.
\bjournal{Numerische Mathematik}
\bvolume{107}
\bpages{181--211}.
\end{barticle}
\endbibitem

\bibitem{Bachoc2013}
\begin{barticle}[author]
\bauthor{\bsnm{Bachoc},~\bfnm{F.}\binits{F.}}
(\byear{2013}).
\btitle{Cross validation and maximum likelihood estimations of hyper-parameters of {G}aussian processes with model misspecification}.
\bjournal{Computational Statistics \& Data Analysis}
\bvolume{66}
\bpages{55--69}.
\end{barticle}
\endbibitem

\bibitem{Bachoc2018}
\begin{barticle}[author]
\bauthor{\bsnm{Bachoc},~\bfnm{F.}\binits{F.}}
(\byear{2018}).
\btitle{Asymptotic analysis of covariance parameter estimation for {G}aussian processes in the misspecified case}.
\bjournal{Bernoulli}
\bvolume{24}
\bpages{1531--1575}.
\end{barticle}
\endbibitem

\bibitem{bachoc2025posterior}
\begin{barticle}[author]
\bauthor{\bsnm{Bachoc},~\bfnm{Fran{\c{c}}ois}\binits{F.}} \AND \bauthor{\bsnm{Lagnoux},~\bfnm{Agn{\`e}s}\binits{A.}}
(\byear{2025}).
\btitle{Posterior contraction rates for constrained deep {Gaussian} processes in density estimation and classification}.
\bjournal{Communications in Statistics-Theory and Methods}
\bvolume{54}
\bpages{774--811}.
\end{barticle}
\endbibitem

\bibitem{Bachoc2019cMLE}
\begin{barticle}[author]
\bauthor{\bsnm{Bachoc},~\bfnm{F.}\binits{F.}}, \bauthor{\bsnm{Lagnoux},~\bfnm{A.}\binits{A.}} \AND \bauthor{\bsnm{L{\'o}pez-Lopera},~\bfnm{A.~F.}\binits{A.~F.}}
(\byear{2019}).
\btitle{Maximum likelihood estimation for {Gaussian} processes under inequality constraints}.
\bjournal{{E}lectronic {J}ournal of {S}tatistics}
\bvolume{13}
\bpages{2921--2969}.
\end{barticle}
\endbibitem

\bibitem{bachoc2022asymptotically}
\begin{barticle}[author]
\bauthor{\bsnm{Bachoc},~\bfnm{Fran{\c{c}}ois}\binits{F.}}, \bauthor{\bsnm{Porcu},~\bfnm{Emilio}\binits{E.}}, \bauthor{\bsnm{Bevilacqua},~\bfnm{Moreno}\binits{M.}}, \bauthor{\bsnm{Furrer},~\bfnm{Reinhard}\binits{R.}} \AND \bauthor{\bsnm{Faouzi},~\bfnm{Tarik}\binits{T.}}
(\byear{2022}).
\btitle{Asymptotically equivalent prediction in multivariate geostatistics}.
\bjournal{Bernoulli}
\bvolume{28}
\bpages{2518--2545}.
\end{barticle}
\endbibitem

\bibitem{BartonPoor1988}
\begin{barticle}[author]
\bauthor{\bsnm{Barton},~\bfnm{R.~J.}\binits{R.~J.}} \AND \bauthor{\bsnm{Poor},~\bfnm{H.~V.}\binits{H.~V.}}
(\byear{1988}).
\btitle{Signal detection in fractional {G}aussian noise}.
\bjournal{IEEE Transactions on Information Theory}
\bvolume{34}
\bpages{943--959}.
\end{barticle}
\endbibitem

\bibitem{Berlinet2004}
\begin{bbook}[author]
\bauthor{\bsnm{Berlinet},~\bfnm{A.}\binits{A.}} \AND \bauthor{\bsnm{Thomas{-}Agnan},~\bfnm{C.}\binits{C.}}
(\byear{2004}).
\btitle{Reproducing Kernel {H}ilbert Spaces in Probability and Statistics}.
\bpublisher{Springer}.
\end{bbook}
\endbibitem

\bibitem{bolin2023equivalence}
\begin{barticle}[author]
\bauthor{\bsnm{Bolin},~\bfnm{David}\binits{D.}} \AND \bauthor{\bsnm{Kirchner},~\bfnm{Kristin}\binits{K.}}
(\byear{2023}).
\btitle{Equivalence of measures and asymptotically optimal linear prediction for {Gaussian} random fields with fractional-order covariance operators}.
\bjournal{Bernoulli}
\bvolume{29}
\bpages{1476--1504}.
\end{barticle}
\endbibitem

\bibitem{castillo2024deep}
\begin{barticle}[author]
\bauthor{\bsnm{Castillo},~\bfnm{Isma{\~A}{\c{G}}l}\binits{I.}} \AND \bauthor{\bsnm{Randrianarisoa},~\bfnm{Thibault}\binits{T.}}
(\byear{2025}+).
\btitle{Deep horseshoe {Gaussian processes}}.
\bjournal{The Annals of Statistics}.
\bnote{To appear}.
\end{barticle}
\endbibitem

\bibitem{Chen2021}
\begin{barticle}[author]
\bauthor{\bsnm{Chen},~\bfnm{Y.}\binits{Y.}}, \bauthor{\bsnm{Owhadi},~\bfnm{H.}\binits{H.}} \AND \bauthor{\bsnm{Stuart},~\bfnm{A.~M.}\binits{A.~M.}}
(\byear{2021}).
\btitle{Consistency of empirical {B}ayes and kernel flow for hierarchical parameter estimation}.
\bjournal{Mathematics of Computation}
\bvolume{90}
\bpages{2527--2578}.
\end{barticle}
\endbibitem

\bibitem{cressie2015statistics}
\begin{bbook}[author]
\bauthor{\bsnm{Cressie},~\bfnm{Noel}\binits{N.}}
(\byear{2015}).
\btitle{Statistics for Spatial Data}.
\bpublisher{John Wiley \& Sons}.
\end{bbook}
\endbibitem

\bibitem{DaVeiga2012GPineqconst}
\begin{barticle}[author]
\bauthor{\bsnm{Da~Veiga},~\bfnm{S\'ebastien}\binits{S.}} \AND \bauthor{\bsnm{Marrel},~\bfnm{Amandine}\binits{A.}}
(\byear{2012}).
\btitle{Gaussian process modeling with inequality constraints}.
\bjournal{Annales de la facult\'e des sciences de Toulouse Math\'ematiques}
\bvolume{21}
\bpages{529--555}.
\end{barticle}
\endbibitem

\bibitem{damianou2013deep}
\begin{binproceedings}[author]
\bauthor{\bsnm{Damianou},~\bfnm{Andreas}\binits{A.}} \AND \bauthor{\bsnm{Lawrence},~\bfnm{Neil~D}\binits{N.~D.}}
(\byear{2013}).
\btitle{Deep {G}aussian processes}.
In \bbooktitle{Proceedings of the Sixteenth International Conference on Artificial Intelligence and Statistics}
\bpages{207--215}.
\bpublisher{PMLR}.
\end{binproceedings}
\endbibitem

\bibitem{DeMarchiSchaback2010}
\begin{barticle}[author]
\bauthor{\bsnm{De~Marchi},~\bfnm{S.}\binits{S.}} \AND \bauthor{\bsnm{Schaback},~\bfnm{R.}\binits{R.}}
(\byear{2010}).
\btitle{Stability of kernel-based interpolation}.
\bjournal{Advances in Computational Mathematics}
\bvolume{32}
\bpages{155--161}.
\end{barticle}
\endbibitem

\bibitem{DickKritzer2022}
\begin{bbook}[author]
\bauthor{\bsnm{Dick},~\bfnm{J.}\binits{J.}}, \bauthor{\bsnm{Kritzer},~\bfnm{P.}\binits{P.}} \AND \bauthor{\bsnm{Pillichshammer},~\bfnm{F.}\binits{F.}}
(\byear{2022}).
\btitle{Lattice Rules: Numerical Integration, Approximation, and Discrepancy}.
\bpublisher{Springer}.
\end{bbook}
\endbibitem

\bibitem{Driscoll1973}
\begin{barticle}[author]
\bauthor{\bsnm{Driscoll},~\bfnm{M.~F.}\binits{M.~F.}}
(\byear{1973}).
\btitle{The reproducing kernel {H}ilbert space structure of the sample paths of a {G}aussian process}.
\bjournal{Zeitschrift f{\"u}r Wahrscheinlichkeitstheorie und Verwandte Gebiete}
\bvolume{26}
\bpages{309--316}.
\end{barticle}
\endbibitem

\bibitem{finocchio2023posterior}
\begin{barticle}[author]
\bauthor{\bsnm{Finocchio},~\bfnm{Gianluca}\binits{G.}} \AND \bauthor{\bsnm{Schmidt-Hieber},~\bfnm{Johannes}\binits{J.}}
(\byear{2023}).
\btitle{Posterior contraction for deep {Gaussian} process priors}.
\bjournal{Journal of Machine Learning Research}
\bvolume{24}
\bpages{1--49}.
\end{barticle}
\endbibitem

\bibitem{genton2001classes}
\begin{barticle}[author]
\bauthor{\bsnm{Genton},~\bfnm{Marc~G}\binits{M.~G.}}
(\byear{2001}).
\btitle{Classes of kernels for machine learning: a statistics perspective}.
\bjournal{Journal of Machine Learning Research}
\bvolume{2}
\bpages{299--312}.
\end{barticle}
\endbibitem

\bibitem{Geoga2023}
\begin{barticle}[author]
\bauthor{\bsnm{Geoga},~\bfnm{C.~J.}\binits{C.~J.}}, \bauthor{\bsnm{Marin},~\bfnm{O.}\binits{O.}}, \bauthor{\bsnm{Schanen},~\bfnm{M.}\binits{M.}} \AND \bauthor{\bsnm{Stein},~\bfnm{M.~L.}\binits{M.~L.}}
(\byear{2023}).
\btitle{Fitting {M}at{\'e}rn smoothness parameters using automatic differentiation}.
\bjournal{Statistics and Computing}
\bvolume{33}
\bpages{48}.
\end{barticle}
\endbibitem

\bibitem{ghosal2017fundamentals}
\begin{bbook}[author]
\bauthor{\bsnm{Ghosal},~\bfnm{Subhashis}\binits{S.}} \AND \bauthor{\bparticle{Van~der} \bsnm{Vaart},~\bfnm{Aad~W}\binits{A.~W.}}
(\byear{2017}).
\btitle{Fundamentals of Nonparametric {B}ayesian Inference}.
\bpublisher{Cambridge University Press}.
\end{bbook}
\endbibitem

\bibitem{Gramacy2020}
\begin{bbook}[author]
\bauthor{\bsnm{Gramacy},~\bfnm{R.~B.}\binits{R.~B.}}
(\byear{2020}).
\btitle{Surrogates: {G}aussian Process Modeling, Design, and Optimization for the Applied Sciences}.
\bpublisher{CRC Press}.
\end{bbook}
\endbibitem

\bibitem{HadjiSzabo2021}
\begin{barticle}[author]
\bauthor{\bsnm{Hadji},~\bfnm{A.}\binits{A.}} \AND \bauthor{\bsnm{Szab\'{o}},~\bfnm{B.}\binits{B.}}
(\byear{2021}).
\btitle{Can we trust {B}ayesian uncertainty quantification from {G}aussian process priors with squared exponential covariance kernel?}
\bjournal{SIAM/ASA Journal on Uncertainty Quantification}
\bvolume{9}
\bpages{185--230}.
\end{barticle}
\endbibitem

\bibitem{Henderson2024}
\begin{barticle}[author]
\bauthor{\bsnm{Henderson},~\bfnm{I.}\binits{I.}}
(\byear{2024}).
\btitle{{S}obolev regularity of {G}aussian random fields}.
\bjournal{Journal of Functional Analysis}
\bvolume{286}
\bpages{110241}.
\end{barticle}
\endbibitem

\bibitem{ibragimov2012gaussian}
\begin{bbook}[author]
\bauthor{\bsnm{Ibragimov},~\bfnm{Ildar~Abdulovich}\binits{I.~A.}} \AND \bauthor{\bsnm{Rozanov},~\bfnm{Yurii~Antol'evich}\binits{Y.~A.}}
(\byear{2012}).
\btitle{{G}aussian Random Processes}.
\bpublisher{Springer}.
\end{bbook}
\endbibitem

\bibitem{Iske2018}
\begin{bbook}[author]
\bauthor{\bsnm{Iske},~\bfnm{A.}\binits{A.}}
(\byear{2018}).
\btitle{Approximation Theory and Algorithms for Data Analysis}.
\bpublisher{Springer}.
\end{bbook}
\endbibitem

\bibitem{Kanagawa2018}
\begin{barticle}[author]
\bauthor{\bsnm{Kanagawa},~\bfnm{M.}\binits{M.}}, \bauthor{\bsnm{Hennig},~\bfnm{P.}\binits{P.}}, \bauthor{\bsnm{Sejdinovic},~\bfnm{D.}\binits{D.}} \AND \bauthor{\bsnm{Sriperumbudur},~\bfnm{B.~K.}\binits{B.~K.}}
(\byear{2018}).
\btitle{Gaussian processes and kernel methods: a review on connections and equivalences}.
\bjournal{arXiv:1807.02582v1}.
\end{barticle}
\endbibitem

\bibitem{Karvonen2023b}
\begin{barticle}[author]
\bauthor{\bsnm{Karvonen},~\bfnm{T.}\binits{T.}}
(\byear{2023}).
\btitle{Small sample spaces for {G}aussian processes}.
\bjournal{Bernoulli}
\bvolume{29}
\bpages{875--900}.
\end{barticle}
\endbibitem

\bibitem{KarvonenOates2023}
\begin{barticle}[author]
\bauthor{\bsnm{Karvonen},~\bfnm{T.}\binits{T.}} \AND \bauthor{\bsnm{Oates},~\bfnm{C.~J.}\binits{C.~J.}}
(\byear{2023}).
\btitle{Maximum likelihood estimation in {G}aussian process regression is ill-posed}.
\bjournal{Journal of Machine Learning Research}
\bvolume{24}
\bpages{1--47}.
\end{barticle}
\endbibitem

\bibitem{KarvonenSantinWenzel2025}
\begin{barticle}[author]
\bauthor{\bsnm{Karvonen},~\bfnm{T.}\binits{T.}}, \bauthor{\bsnm{Santin},~\bfnm{G.}\binits{G.}} \AND \bauthor{\bsnm{Wenzel},~\bfnm{T.}\binits{T.}}
(\byear{2025}).
\btitle{General superconvergence for kernel-based approximation}.
\bjournal{arXiv:2505.11435v1}.
\end{barticle}
\endbibitem

\bibitem{Karvonen2020}
\begin{barticle}[author]
\bauthor{\bsnm{Karvonen},~\bfnm{T.}\binits{T.}}, \bauthor{\bsnm{Wynne},~\bfnm{G.}\binits{G.}}, \bauthor{\bsnm{Tronarp},~\bfnm{F.}\binits{F.}}, \bauthor{\bsnm{Oates},~\bfnm{C.~J.}\binits{C.~J.}} \AND \bauthor{\bsnm{S{\"a}rkk{\"a}},~\bfnm{S.}\binits{S.}}
(\byear{2020}).
\btitle{Maximum likelihood estimation and uncertainty quantification for {G}aussian process approximation of deterministic functions}.
\bjournal{SIAM/ASA Journal on Uncertainty Quantification}
\bvolume{8}
\bpages{926--958}.
\end{barticle}
\endbibitem

\bibitem{KorteStapff2024}
\begin{barticle}[author]
\bauthor{\bsnm{Korte-Stapff},~\bfnm{M.}\binits{M.}}, \bauthor{\bsnm{Karvonen},~\bfnm{T.}\binits{T.}} \AND \bauthor{\bsnm{Moulines},~\bfnm{E.}\binits{E.}}
(\byear{2025}).
\btitle{Smoothness estimation for {W}hittle-{M}at\'{e}rn processes on closed {R}iemannian manifolds}.
\bjournal{Stochastic Processes and Their Applications}
\bvolume{189}.
\bnote{Article no. 104685}.
\end{barticle}
\endbibitem

\bibitem{LukicBeder2001}
\begin{barticle}[author]
\bauthor{\bsnm{Luki{\'c}},~\bfnm{M.~N.}\binits{M.~N.}} \AND \bauthor{\bsnm{Beder},~\bfnm{J.~H.}\binits{J.~H.}}
(\byear{2001}).
\btitle{Stochastic processes with sample paths in reproducing kernel {H}ilbert spaces}.
\bjournal{Transactions of the American Mathematical Society}
\bvolume{353}
\bpages{3945--3969}.
\end{barticle}
\endbibitem

\bibitem{Narcowich2006}
\begin{barticle}[author]
\bauthor{\bsnm{Narcowich},~\bfnm{F.~J.}\binits{F.~J.}}, \bauthor{\bsnm{Ward},~\bfnm{J.~D.}\binits{J.~D.}} \AND \bauthor{\bsnm{Wendland},~\bfnm{H.}\binits{H.}}
(\byear{2006}).
\btitle{{S}obolev error estimates and a {B}ernstein inequality for scattered data interpolation via radial basis functions}.
\bjournal{Constructive Approximation}
\bvolume{24}
\bpages{175--186}.
\end{barticle}
\endbibitem

\bibitem{Naslidnyk2024}
\begin{barticle}[author]
\bauthor{\bsnm{Naslidnyk},~\bfnm{M.}\binits{M.}}, \bauthor{\bsnm{Kanagawa},~\bfnm{M.}\binits{M.}}, \bauthor{\bsnm{Karvonen},~\bfnm{T.}\binits{T.}} \AND \bauthor{\bsnm{Mahsereci},~\bfnm{M.}\binits{M.}}
(\byear{2025}).
\btitle{Comparing scale parameter estimators for {G}aussian process interpolation with the {B}rownian motion prior: Leave-one-out cross validation and maximum likelihood}.
\bjournal{SIAM/ASA Journal on Uncertainty Quantification}
\bvolume{13}
\bpages{679--717}.
\end{barticle}
\endbibitem

\bibitem{nickisch2008approximations}
\begin{barticle}[author]
\bauthor{\bsnm{Nickisch},~\bfnm{Hannes}\binits{H.}}, \bauthor{\bsnm{Rasmussen},~\bfnm{Carl~Edward}\binits{C.~E.}} \betal{et~al.}
(\byear{2008}).
\btitle{Approximations for binary {Gaussian} process classification}.
\bjournal{Journal of Machine Learning Research}
\bvolume{9}
\bpages{2035--2078}.
\end{barticle}
\endbibitem

\bibitem{Niederreiter1992}
\begin{bbook}[author]
\bauthor{\bsnm{Niederreiter},~\bfnm{H.}\binits{H.}}
(\byear{1992}).
\btitle{Random Number Generation and Quasi-{M}onte {C}arlo Methods}.
\bpublisher{Society for Industrial and Applied Mathematics}.
\end{bbook}
\endbibitem

\bibitem{Novak1988}
\begin{bbook}[author]
\bauthor{\bsnm{Novak},~\bfnm{E.}\binits{E.}}
(\byear{1988}).
\btitle{Deterministic and Stochastic Error Bounds in Numerical Analysis}.
\bpublisher{Springer-Verlag}.
\end{bbook}
\endbibitem

\bibitem{NovakWozniakowski2008}
\begin{bbook}[author]
\bauthor{\bsnm{Novak},~\bfnm{E.}\binits{E.}} \AND \bauthor{\bsnm{Wo{\'z}niakowski},~\bfnm{H.}\binits{H.}}
(\byear{2008}).
\btitle{Tractability of Multivariate Problems. Volume {I}: Linear Information}.
\bpublisher{European Mathematical Society}.
\end{bbook}
\endbibitem

\bibitem{NovakWozniakowski2010}
\begin{bbook}[author]
\bauthor{\bsnm{Novak},~\bfnm{E.}\binits{E.}} \AND \bauthor{\bsnm{Wo{\'z}niakowski},~\bfnm{H.}\binits{H.}}
(\byear{2010}).
\btitle{Tractability of Multivariate Problems. Volume {II}: Standard Information for Functionals}.
\bpublisher{European Mathematical Society}.
\end{bbook}
\endbibitem

\bibitem{Paulsen2016}
\begin{bbook}[author]
\bauthor{\bsnm{Paulsen},~\bfnm{V.~I.}\binits{V.~I.}} \AND \bauthor{\bsnm{Raghupathi},~\bfnm{M.}\binits{M.}}
(\byear{2016}).
\btitle{An Introduction to the Theory of Reproducing Kernel {H}ilbert Spaces}.
\bpublisher{Cambridge University Press}.
\end{bbook}
\endbibitem

\bibitem{Petit2023}
\begin{barticle}[author]
\bauthor{\bsnm{Petit},~\bfnm{S.}\binits{S.}}
(\byear{2025}).
\btitle{An asymptotic study of the joint maximum likelihood estimation of the regularity and the amplitude parameters of a periodized {M}atérn model}.
\bjournal{Electronic Journal of Statistics}
\bvolume{19}
\bpages{2052--2094}.
\end{barticle}
\endbibitem

\bibitem{PronzatoZhigljavsky2023}
\begin{barticle}[author]
\bauthor{\bsnm{Pronzato},~\bfnm{L.}\binits{L.}} \AND \bauthor{\bsnm{Zhigljavsky},~\bfnm{A.}\binits{A.}}
(\byear{2023}).
\btitle{Quasi-uniform designs with optimal and near-optimal uniformity constant}.
\bjournal{Journal of Approximation Theory}
\bvolume{294}.
\bnote{Article no.\ 105931}.
\end{barticle}
\endbibitem

\bibitem{putter2001effect}
\begin{barticle}[author]
\bauthor{\bsnm{Putter},~\bfnm{Hein}\binits{H.}} \AND \bauthor{\bsnm{Young},~\bfnm{G.~A.}\binits{G.~A.}}
(\byear{2001}).
\btitle{On the effect of covariance function estimation on the accuracy of kriging predictors}.
\bjournal{Bernoulli}
\bvolume{7}
\bpages{421--438}.
\end{barticle}
\endbibitem

\bibitem{RasmussenWilliams2006}
\begin{bbook}[author]
\bauthor{\bsnm{Rasmussen},~\bfnm{C.~E.}\binits{C.~E.}} \AND \bauthor{\bsnm{Williams},~\bfnm{C.~K.~I.}\binits{C.~K.~I.}}
(\byear{2006}).
\btitle{Gaussian Processes for Machine Learning}.
\bpublisher{MIT Press}.
\end{bbook}
\endbibitem

\bibitem{Pallavi2019BayesianShapeGPs}
\begin{barticle}[author]
\bauthor{\bsnm{{Ray}},~\bfnm{Pallavi}\binits{P.}}, \bauthor{\bsnm{{Pati}},~\bfnm{Debdeep}\binits{D.}} \AND \bauthor{\bsnm{{Bhattacharya}},~\bfnm{Anirban}\binits{A.}}
(\byear{2020}).
\btitle{Efficient {Bayesian} shape-restricted function estimation with constrained {Gaussian} process priors}.
\bjournal{Statistics and Computing}
\bvolume{30}
\bpages{839-853}.
\end{barticle}
\endbibitem

\bibitem{Sacks1989}
\begin{barticle}[author]
\bauthor{\bsnm{Sacks},~\bfnm{J.}\binits{J.}}, \bauthor{\bsnm{Welch},~\bfnm{W.~J.}\binits{W.~J.}}, \bauthor{\bsnm{Mitchell},~\bfnm{T.~J.}\binits{T.~J.}} \AND \bauthor{\bsnm{Wynn},~\bfnm{H.~P.}\binits{H.~P.}}
(\byear{1989}).
\btitle{Design and analysis of computer experiments}.
\bjournal{Statistical Science}
\bvolume{4}
\bpages{409--435}.
\end{barticle}
\endbibitem

\bibitem{Santner2003}
\begin{bbook}[author]
\bauthor{\bsnm{Santner},~\bfnm{T.~J.}\binits{T.~J.}}, \bauthor{\bsnm{Williams},~\bfnm{B.~J.}\binits{B.~J.}} \AND \bauthor{\bsnm{Notz},~\bfnm{W.~I.}\binits{W.~I.}}
(\byear{2003}).
\btitle{The Design and Analysis of Computer Experiments}.
\bpublisher{Springer}.
\end{bbook}
\endbibitem

\bibitem{SanzAlonso2024}
\begin{barticle}[author]
\bauthor{\bsnm{Sanz-Alonso},~\bfnm{D.}\binits{D.}} \AND \bauthor{\bsnm{Yang},~\bfnm{R.}\binits{R.}}
(\byear{2025}).
\btitle{{G}aussian process regression under computational and epistemic misspecification}.
\bjournal{SIAM Journal on Numerical Analysis}
\bvolume{63}
\bpages{495--519}.
\end{barticle}
\endbibitem

\bibitem{SchabackWendland2001}
\begin{bincollection}[author]
\bauthor{\bsnm{Schaback},~\bfnm{R.}\binits{R.}} \AND \bauthor{\bsnm{Wendland},~\bfnm{H.}\binits{H.}}
(\byear{2001}).
\btitle{Approximation by positive definite kernels}.
In \bbooktitle{Advanced Problems in Constructive Approximation}
\bpages{203--222}.
\end{bincollection}
\endbibitem

\bibitem{Scheuerer2010}
\begin{barticle}[author]
\bauthor{\bsnm{Scheuerer},~\bfnm{M.}\binits{M.}}
(\byear{2010}).
\btitle{Regularity of the sample paths of a general second order random field}.
\bjournal{Stochastic Processes and Their Applications}
\bvolume{120}
\bpages{1879--1897}.
\end{barticle}
\endbibitem

\bibitem{skorokhod1973absolute}
\begin{barticle}[author]
\bauthor{\bsnm{Skorokhod},~\bfnm{Anatolii~Volodimirovich}\binits{A.~V.}} \AND \bauthor{\bsnm{Yadrenko},~\bfnm{Mikha{\i}lo~Iosipovich}\binits{M.~I.}}
(\byear{1973}).
\btitle{On absolute continuity of measures corresponding to homogeneous {Gaussian} fields}.
\bjournal{Theory of Probability \& Its Applications}
\bvolume{18}
\bpages{27--40}.
\end{barticle}
\endbibitem

\bibitem{SniekersVaart2015a}
\begin{barticle}[author]
\bauthor{\bsnm{Sniekers},~\bfnm{S.}\binits{S.}} \AND \bauthor{\bparticle{van~der} \bsnm{Vaart},~\bfnm{A.}\binits{A.}}
(\byear{2015}).
\btitle{Adaptive {B}ayesian credible sets in regression with a {G}aussian process prior}.
\bjournal{Electronic Journal of Statistics}
\bvolume{9}
\bpages{2475--2527}.
\end{barticle}
\endbibitem

\bibitem{SniekersVaart2015b}
\begin{barticle}[author]
\bauthor{\bsnm{Sniekers},~\bfnm{S.}\binits{S.}} \AND \bauthor{\bparticle{van~der} \bsnm{Vaart},~\bfnm{A.}\binits{A.}}
(\byear{2015}).
\btitle{Credible sets in the fixed design model with {B}rownian motion prior}.
\bjournal{Journal of Statistical Planning and Inference}
\bvolume{166}
\bpages{78--86}.
\end{barticle}
\endbibitem

\bibitem{stein1990uniform}
\begin{barticle}[author]
\bauthor{\bsnm{Stein},~\bfnm{Michael}\binits{M.}}
(\byear{1990}).
\btitle{Uniform asymptotic optimality of linear predictions of a random field using an incorrect second-order structure}.
\bjournal{The Annals of Statistics}
\bpages{850--872}.
\end{barticle}
\endbibitem

\bibitem{stein1988asymptotically}
\begin{barticle}[author]
\bauthor{\bsnm{Stein},~\bfnm{Michael~L}\binits{M.~L.}}
(\byear{1988}).
\btitle{Asymptotically efficient prediction of a random field with a misspecified covariance function}.
\bjournal{The Annals of Statistics}
\bpages{55--63}.
\end{barticle}
\endbibitem

\bibitem{stein1990bounds}
\begin{barticle}[author]
\bauthor{\bsnm{Stein},~\bfnm{Michael~L}\binits{M.~L.}}
(\byear{1990}).
\btitle{Bounds on the efficiency of linear predictions using an incorrect covariance function}.
\bjournal{The Annals of Statistics}
\bpages{1116--1138}.
\end{barticle}
\endbibitem

\bibitem{Stein1990}
\begin{barticle}[author]
\bauthor{\bsnm{Stein},~\bfnm{M.~L.}\binits{M.~L.}}
(\byear{1990}).
\btitle{A comparison of generalized cross validation and modified maximum likelihood for estimating the parameters of a stochastic process}.
\bjournal{The Annals of Statistics}
\bvolume{18}
\bpages{1139--1157}.
\end{barticle}
\endbibitem

\bibitem{Stein1999}
\begin{bbook}[author]
\bauthor{\bsnm{Stein},~\bfnm{M.~L.}\binits{M.~L.}}
(\byear{1999}).
\btitle{Interpolation of Spatial Data: {S}ome Theory for Kriging}.
\bpublisher{Springer}.
\end{bbook}
\endbibitem

\bibitem{Steinwart2019}
\begin{barticle}[author]
\bauthor{\bsnm{Steinwart},~\bfnm{I.}\binits{I.}}
(\byear{2019}).
\btitle{Convergence types and rates in generic {K}arhunen-{L}oève expansions with applications to sample path properties}.
\bjournal{Potential Analysis}
\bvolume{51}
\bpages{361--395}.
\end{barticle}
\endbibitem

\bibitem{SteinwartScovel2012}
\begin{barticle}[author]
\bauthor{\bsnm{Steinwart},~\bfnm{I.}\binits{I.}} \AND \bauthor{\bsnm{Scovel},~\bfnm{C.}\binits{C.}}
(\byear{2012}).
\btitle{{M}ercer's theorem on general domains: {O}n the interaction between measures, kernels, and {RKHS}s}.
\bjournal{Constructive Approximation}
\bvolume{35}
\bpages{363--417}.
\end{barticle}
\endbibitem

\bibitem{szabo2025adaptation}
\begin{barticle}[author]
\bauthor{\bsnm{Szabo},~\bfnm{Botond}\binits{B.}}, \bauthor{\bsnm{Hadji},~\bfnm{Amine}\binits{A.}} \AND \bauthor{\bparticle{van~der} \bsnm{Vaart},~\bfnm{Aad}\binits{A.}}
(\byear{2025}).
\btitle{Adaptation using spatially distributed {Gaussian} processes}.
\bjournal{Journal of the American Statistical Association}
\bpages{1--34}.
\end{barticle}
\endbibitem

\bibitem{szabo2015frequentist}
\begin{barticle}[author]
\bauthor{\bsnm{Szab{\'o}},~\bfnm{BT}\binits{B.}}, \bauthor{\bparticle{van~der} \bsnm{Vaart},~\bfnm{AW}\binits{A.}} \AND \bauthor{\bparticle{van} \bsnm{Zanten},~\bfnm{JH}\binits{J.}}
(\byear{2015}).
\btitle{Frequentist coverage of adaptive nonparametric {Bayesian} credible sets}.
\bjournal{The Annals of Statistics}
\bvolume{43}
\bpages{1391--1428}.
\end{barticle}
\endbibitem

\bibitem{Teckentrup2020}
\begin{barticle}[author]
\bauthor{\bsnm{Teckentrup},~\bfnm{A.~L.}\binits{A.~L.}}
(\byear{2020}).
\btitle{Convergence of {G}aussian process regression with estimated hyper-parameters and applications in {B}ayesian inverse problems}.
\bjournal{SIAM/ASA Journal on Uncertainty Quantification}
\bvolume{8}
\bpages{1310--1337}.
\end{barticle}
\endbibitem

\bibitem{VaartZanten2011}
\begin{barticle}[author]
\bauthor{\bparticle{van~der} \bsnm{Vaart},~\bfnm{A.}\binits{A.}} \AND \bauthor{\bparticle{van} \bsnm{Zanten},~\bfnm{H.}\binits{H.}}
(\byear{2011}).
\btitle{Information rates of nonparametric {G}aussian process methods}.
\bjournal{Journal of Machine Learning Research}
\bvolume{12}
\bpages{2095--2119}.
\end{barticle}
\endbibitem

\bibitem{van2008rates}
\begin{barticle}[author]
\bauthor{\bparticle{van~der} \bsnm{Vaart},~\bfnm{AW}\binits{A.}} \AND \bauthor{\bparticle{van} \bsnm{Zanten},~\bfnm{JH}\binits{J.}}
(\byear{2008}).
\btitle{Rates of contraction of posterior distributions based on {Gaussian} process priors}.
\bjournal{The Annals of Statistics}
\bvolume{36}
\bpages{1435--1463}.
\end{barticle}
\endbibitem

\bibitem{van2009adaptive}
\begin{barticle}[author]
\bauthor{\bparticle{van~der} \bsnm{Vaart},~\bfnm{AW}\binits{A.}} \AND \bauthor{\bparticle{van} \bsnm{Zanten},~\bfnm{JH}\binits{J.}}
(\byear{2009}).
\btitle{Adaptive Bayesian estimation using a {Gaussian} random field with inverse {Gamma} bandwitdh}.
\bjournal{The Annals of Statistics}
\bvolume{37}
\bpages{2655--2675}.
\end{barticle}
\endbibitem

\bibitem{VaartZanten2008}
\begin{binbook}[author]
\bauthor{\bparticle{van~der} \bsnm{Vaart},~\bfnm{A.~W.}\binits{A.~W.}} \AND \bauthor{\bparticle{van} \bsnm{Zanten},~\bfnm{J.~H.}\binits{J.~H.}}
(\byear{2008}).
\btitle{Reproducing Kernel {H}ilbert spaces of {G}aussian Priors}.
In \bbooktitle{Pushing the Limits of Contemporary Statistics: Contributions in Honor of Jayanta K.\ Ghosh}.
\bseries{IMS Collections}
\bvolume{3}
\bpages{200--222}.
\bpublisher{Institute of Mathematical Statistics}.
\end{binbook}
\endbibitem

\bibitem{Wahba1990}
\begin{bbook}[author]
\bauthor{\bsnm{Wahba},~\bfnm{G.}\binits{G.}}
(\byear{1990}).
\btitle{Spline Models for Observational Data}.
\bpublisher{Society for Industrial and Applied Mathematics}.
\end{bbook}
\endbibitem

\bibitem{Wang2019DiffEq}
\begin{barticle}[author]
\bauthor{\bsnm{Wang},~\bfnm{Junyang}\binits{J.}}, \bauthor{\bsnm{Cockayne},~\bfnm{Jon}\binits{J.}} \AND \bauthor{\bsnm{Oates},~\bfnm{Chris.~J.}\binits{C.~J.}}
(\byear{2020}).
\btitle{A role for symmetry in the {B}ayesian solution of differential equations}.
\bjournal{Bayesian Analysis}
\bvolume{15}
\bpages{1057--1085}.
\end{barticle}
\endbibitem

\bibitem{WangJing2022}
\begin{barticle}[author]
\bauthor{\bsnm{Wang},~\bfnm{W.}\binits{W.}} \AND \bauthor{\bsnm{Jing},~\bfnm{B.~Y.}\binits{B.~Y.}}
(\byear{2022}).
\btitle{{G}aussian process regression: Optimality, robustness, and relationship with kernel ridge regression}.
\bjournal{Journal of Machine Learning Research}
\bvolume{23}
\bpages{1--67}.
\end{barticle}
\endbibitem

\bibitem{Wendland2005}
\begin{bbook}[author]
\bauthor{\bsnm{Wendland},~\bfnm{H.}\binits{H.}}
(\byear{2005}).
\btitle{Scattered Data Approximation}.
\bpublisher{Cambridge University Press}.
\end{bbook}
\endbibitem

\bibitem{Wynne2021}
\begin{barticle}[author]
\bauthor{\bsnm{Wynne},~\bfnm{G.}\binits{G.}}, \bauthor{\bsnm{Briol},~\bfnm{F.~X.}\binits{F.~X.}} \AND \bauthor{\bsnm{Girolami},~\bfnm{M.}\binits{M.}}
(\byear{2021}).
\btitle{Convergence guarantees for {G}aussian process means with misspecified likelihoods and smoothness}.
\bjournal{Journal of Machine Learning Research}
\bvolume{22}
\bpages{1--40}.
\end{barticle}
\endbibitem

\bibitem{xu2017tukey}
\begin{barticle}[author]
\bauthor{\bsnm{Xu},~\bfnm{Ganggang}\binits{G.}} \AND \bauthor{\bsnm{Genton},~\bfnm{Marc~G}\binits{M.~G.}}
(\byear{2017}).
\btitle{{T}ukey $g$-and-$h$ random fields}.
\bjournal{Journal of the American Statistical Association}
\bvolume{112}
\bpages{1236--1249}.
\end{barticle}
\endbibitem

\bibitem{XuStein2017}
\begin{barticle}[author]
\bauthor{\bsnm{Xu},~\bfnm{W.}\binits{W.}} \AND \bauthor{\bsnm{Stein},~\bfnm{M.~L.}\binits{M.~L.}}
(\byear{2017}).
\btitle{Maximum likelihood estimation for a smooth {G}aussian random field model}.
\bjournal{SIAM/ASA Journal on Uncertainty Quantification}
\bvolume{5}
\bpages{138--175}.
\end{barticle}
\endbibitem

\bibitem{zhang2004inconsistent}
\begin{barticle}[author]
\bauthor{\bsnm{Zhang},~\bfnm{Hao}\binits{H.}}
(\byear{2004}).
\btitle{Inconsistent estimation and asymptotically equal interpolations in model-based geostatistics}.
\bjournal{Journal of the American Statistical Association}
\bvolume{99}
\bpages{250--261}.
\end{barticle}
\endbibitem

\bibitem{zhang2015doesn}
\begin{barticle}[author]
\bauthor{\bsnm{Zhang},~\bfnm{Hao}\binits{H.}} \AND \bauthor{\bsnm{Cai},~\bfnm{Wenxiang}\binits{W.}}
(\byear{2015}).
\btitle{When doesn't cokriging outperform kriging?}
\bjournal{Statistical Science}
\bpages{176--180}.
\end{barticle}
\endbibitem

\end{thebibliography}
\end{document}